\documentclass{article}

\usepackage{amsmath}
\usepackage{amsfonts}
\usepackage{amssymb}
\usepackage{amsthm}
\usepackage{stmaryrd}
\usepackage[T1]{fontenc}
\usepackage{enumerate}
\usepackage{graphics}
\usepackage{subfig,cases}
\usepackage{empheq}
\usepackage{multirow,booktabs}
\usepackage{float}

\usepackage{tikz}
\usepackage[linesnumbered,boxed, ruled]{algorithm2e}
\SetKwRepeat{Do}{do}{while}
\usepackage{multicol}
\usepackage[colorlinks]{hyperref}

\newcommand{\ds}{\displaystyle}
\newcommand{\pa}{\partial}

\newcommand{\keff}{{  \rm k_{\rm eff}}} 
\newcommand{\keffnp}{{  \rm k_{\rm eff}^{(n+1)}}}
\newcommand{\keffn}{{  \rm k_{\rm eff}^{(n)}}}

\newcommand{\keffref}{{  \rm k_{\rm eff}^{\rm ref}}}

\DeclareMathOperator{\grad}{\bf grad}
\DeclareMathOperator{\dive}{\rm div}

\def\nvec{{\bf n}}

\def\pvec{{\bf p}}
\def\qvec{{\bf q}}

\def\xvec{{\bf x}}

\def\Hvec{{\bf H}}

\def\Lvec{{\bf L}}
\def\xvec{{\bf x}}
\def\Qvec{{\bf Q}}

\def\Bcal{{\bf \mathcal B}}

\def\Ecal{{\bf \mathcal E}}

\def\Ical{{\bf \mathcal I}}
\def\Lcal{{\bf \mathcal L}}

\def\Pcal{{\bf \mathcal P}}

\def\Tcal{{\bf \mathcal T}}

\def\Jcal{{\bf \mathcal J}}
\def\Dbb{{\mathbb D}}
\def\Mbb{{\mathbb M}}

\def\T{\mathbb{T}}
\def\D{\mathbb{D}}
\def\Wud{{\underline W}}
\newcommand{\N}{\mathbb{N}}
\newcommand{\R}{\mathbb{R}}

\newtheorem{remark}{Remark}

\usepackage{authblk}
\begin{document}




\title{On Physics-Based Loss Scaling for MF-PINNs applied to the neutron diffusion equation} 

%
\author[1]{Minh Hieu Do}
\author[1]{Fran\c{c}ois Madiot}
\author[1]{Karim Ammar}
\author[1]{Nicolas G\'erard Castaing}
\affil[1]{Universit\'e Paris-Saclay, CEA, Service d'\'Etudes des R\'eacteurs et de Math\'ematiques Appliqu\'ees, 91191, Gif-sur-Yvette, France.}

\maketitle
\begin{abstract}
Physics-Based Loss Scaling (PBLS) is introduced  for Mixed-Formulation PINNs (MF-PINNs) applied to the  neutron diffusion equation. In particular, we propose a new \textit{scaled} loss function based on the material cross sections, which is equivalent to the classical MF-PINN loss, but  accelerates the convergence and improves accuracy of MF-PINNs. Several numerical experiments on both the fixed source and the k-eigenvalue problem, from one-group to multigroup cases and from two-dimensional (2D) to three-dimensional (3D) configurations, illustrate the efficiency of the proposed scaling method.

\end{abstract}



%
%




\section{Introduction}
Traditional deep learning methods for solving partial differential equations (PDEs) usually reply on purely data-driven models where neural networks are trained in  a supervised learning framework with large amounts of labeled datasets from high-fidelity numerical simulations with traditional numerical methods, experiments, or measurements from sensors. However, this strategy suffers from several limitations. For example, it can ignore the physical  properties such as gorverning equations or conservation laws and leads to  unphysical predictions. Data-driven approaches often struggle with noises in data, which is very common in real world and engineering applications. To overcome these challenges,
Physics-Informed Neural Networks (PINNs), originally introduced  in \cite{raissi2019physics}
by embedding governing equations directly into the loss function  through automatic differentiation of deep neural networks, provide a flexible, mesh-free, and differentiable framework  for the numerical approximation of PDEs while respecting the  physical constraints.
Therefore, PINNs provide a powerful framework for computational science across various domains \cite{cai2021physics,cuomo2022scientific,Mishra2024}. Several variants of PINNs have been investigated in recent years, such as VPINNS \cite{kharazmi2019variational} based on the variational formulation, wPINNs \cite{de2024wpinns} for entropy solutions of hyperbolic
conservation laws, and fPINNs \cite{pang2019fpinns} for space-time fractional
advection-diffusion equations.

Reactor core analysis, which typically relies on solving the neutron transport equation or  the neutron diffusion equation (NDE),  is essential to ensure the safe and reliable operation of a nuclear reactor. Although PINNs have been successfully applied to solve the neutron transport equation in \cite{xie2024boundary, geng2025physics}, the present work focuses on their application to  the neutron diffusion equation. The NDE is a widely used approximation that significantly reduces computational cost while maintaining sufficient accuracy for predicting neutron flux distributions and reactivity at the reactor core scale. In fact, several studies have investigated the application of PINNs to the neutron diffusion equation (NDE) in its primal formulation
 \cite{elhareef2023physics,yang2023physics, bi2025fc,lee2026development,yu2026solving}. It is important to note that the neutron deterministic calculation is normally performed using a two-step scheme \cite{Hebert2020,viallon2024numerical}: a two-dimensional (2D) lattice calculation to generate the homogenized cross sections and a three-dimensional (3D) full-core calculation using the neutron diffusion equation (NDE) with these homogenized cross sections to produce the power distribution. Hence, solving the neutron diffusion equation with heterogeneous media, where homogenized cross sections vary across different regions, is typically required for accurate reactor core modeling.
In this case, the PDE residual loss function of the neutron diffusion equation in primal form is not well defined in the whole domain. Therefore, the interface conditions are incorporated  into the loss function of PINNs to ensure the continuity at the interfaces, which is similar to the Discontinuous Least-squrares Finite Element Methods (DLSFEMs). However, in nuclear
reactor core simulations, the presence of multiple materials  (e.g., fuel, coolant, moderator) leads to several interfaces at the interactions among materials. On the one hand, perfectly training multiple interfaces is challenging, especially in the multigroup case. On the other hand, including interface terms increases the computational time for the total loss function, since the calculation of gradient flux  via backpropagation is required for each interface term. More importantly, the neutron diffusion equation in heterogeneous media frequently often exhibits low-regularity solutions, making it very challenging for primal-form PINNs, which usually perform best with smooth solutions. To address these challenges, PINNs based on the mixed dual form of the one-group neutron diffusion equation is proposed in \cite{do2025physics} where a neural network is trained to produce both the neutron flux and the neutron current. This framework is then extended to the more general setting of the multigroup neutron diffusion and 
is referred to as Mixed Formulation PINNs (MF-PINNs) in \cite{do2026MF}. Compared to the primal form PINNs, MF-PINNs offer several advantages. First of all, the loss function of MF-PINNs is well defined over the entire domain, eliminating the need for interface terms. Secondly, imposing hard boundary conditions to ensure the neural network totally satisfies the boundary condition is much easier since  both flux and current are provided by  MF-PINNs. Moreover, the mixed formulation converts the problem into a first-order system, which reduces its training cost compared to higher-order derivatives of the primal form. Several studies have already explored this direction for other applications that confirm the effectiveness of this method, including the mixed residual method (MIM) \cite{lyu2022mim}, Fourier-based mixed physics informed neural networks (FMPINN) \cite{wu2024solving}
and deep first order system least squares method (FOSLS)
\cite{cai2020deep,bersetche2023deep}.  Furthermore, error analysis for deep mixed residual method is also investigated in \cite{li2024priori} and \cite{gu2024error}. 

It is important to point out that MF-PINNs introduce more residual terms, particularly for the multigroup neutron diffusion equation, which leads to the issue of imbalanced gradient among the multiple groups or terms.
In fact, PINNs minimize a composite (multi-term) objective function that combines several contributions, including interior PDE residuals  related to different physical equations, boundary conditions, interface constraints, data-fitting terms, auxiliary equations, and additional terms.
MF-PINNs involve several coupled equations that are enforced simultaneously, which can produce residual terms with significantly different scales and physical units, especially in the presence of heterogeneous coefficients. Without an appropriate normalization or weighting method, this gradient imbalance can lead to poor convergence, as the optimization process can become dominated by certain groups or terms, preventing the neural network from accurately satisfying all physical constraints.

Various strategies have been proposed to mitigate the loss imbalance in PINNs. In particular, several adaptive weighting approaches have been developed, such as: learning rate annealing (LRA) \cite{wang2021understanding}, neural tangent kernel (NTK) \cite{wang2022NTK}, Augmented Lagrangian relaxation (AL) \cite{son2023AL}, residual-based attention (RBA) \cite{anagnostopoulos2024RBA} and the balanced residual decay rate (BRDR) \cite{chen2025BRDR}.
Rather than proposing  another adaptive weighting scheme for MF-PINNs. Inspired by a  posteriori error estimates for the mixed finite element method of the multigroup neutron diffusion equation \cite{ciarletChapter5}, we would like to propose the Physics Based Loss Scaling (PBLS) method, which scales the standard loss function of MF-PINNs by using only material cross sections without requiring any additional training or complex weighting strategies.
This proper scaling is very simple but significantly  improves physical consistency and solution accuracy in heterogeneous multigroup diffusion problems. MF-PINNs exhibit faster convergence when trained with the scaled loss function.

The rest of the paper is structured as follows. Section \ref{section-NDE} recall the multigroup neutron diffusion equation in both the primal and mixed form. Next,
the MF-PINNs for the source problem is explained in Section \ref{section-pinns-SP}. More importantly, in this section, we introduce Physics-Based Loss Scaling (PBLS) for MF-PINNs. Next, MF-PINNs with scaled loss for  the k-eigenvalue problem 
are discussed in Section \ref{section-MF-PINNs-EP}. Section \ref{section:numerical} is dedicated to numerical test cases for both the fixed source and the k-eigenvalue problem. Finally, some conclusions and perspectives are discussed in Section \ref{section-conclusion}.

The implementation of the proposed scaling method for MF-PINNs  applied to the multigroup neutron diffusion equations is  available at:
\begin{center}
\url{https://github.com/ML-SERMA/PINNs.git}
\end{center}

\section{Neutron diffusion equations }
\label{section-NDE}
The steady state multigroup neutron diffusion equations, which are obtained from the neutron transport equation under the diffusion approximation and multigroup energy discretization,  are commonly used in nuclear reactor neutronics calculations. Let $\Omega$ be a bounded, connected, and open subset of $\mathbb{R}^d$ with $d=2,3$, having a Lipschitz boundary that is piecewise smooth. The boundary condition is now split into three disjoint open parts such that $\partial \Omega = \overline{\Gamma}_{D}\cup \overline{\Gamma}_{N} \cup \overline{\Gamma}_{R} $.  Let $\nvec$ be the unit outward normal vector field to $\pa {\Omega}$. For $G\ge 2$, let $\Ical_G := \{1,...,G\}$ be the set of energy groups. In general, the multigroup neutron diffusion equations can be written as follows, for all $g\in\Ical_G$:

\begin{equation}
\label{eq:primal-EP}
\left\lbrace
\begin{aligned}
    -\dive( D^g \grad \phi^g) +\Sigma_r^g \phi^g -\sum_{g'\neq g}\Sigma_s^{g'\to g} \phi^{g'} &= \frac{\chi^g}{\keff} \sum_{g'=1}^{G}\nu\Sigma_f^{g'} \phi^{g'} \quad  \text{ in }\Omega, \\
    \phi^g &= 0 \quad\text{on} \quad \Gamma_D,\\
    D^g \grad \phi^g\cdot\nvec &= 0  \quad \text{on} \quad\Gamma_N,\\
    D^g \grad \phi^g\cdot\nvec + \frac{1}{2} \phi^g &=0  \quad \text{on} \quad\Gamma_R,\\
\end{aligned}
\right.
\end{equation}

where 
\begin{itemize}
    \item $D^g$  is the  diffusion coefficient of the energy group $g$
    \item $\Sigma_t^g$ is the macroscopic total cross section of the energy group $g$.
    \item $\Sigma_s^{g'\to g}$ is the macroscopic scattering cross sections from the energy group $g'$ to the energy group $g$
    \item $\Sigma_r^g = \Sigma_t^g-\Sigma_s^{g\to g}$ is the macroscopic removal cross section of the energy group $g$
    \item $\chi^g$ is the total fission spectrum of the energy group $g$
    \item $\nu^{g'}$ is the average number of neutrons emitted by fission in the energy group $g'$ and $\Sigma_f^{g'}$ is the macroscopic fission cross section of the energy group $g'$
    \item $\phi^g$ stands for the neutron scalar flux of the energy group $g$
    and  $\keff$ is the multiplication factor
\end{itemize}
 The physical state of the core reactor
is characterized by the effective multiplication factor: if $k_{\text{eff}}=1$: the nuclear chain
reaction is self-sustaining; if $k_{\text{eff}} > 1$: the chain reaction is diverging; if $k_{\text{eff}}<1$: the chain reaction
vanishes.

Let $\T_e$ be the even removal matrix. At (almost) every point in $\Omega$, it is  a matrix of $\R^{G\times G}$:
\begin{equation*}
 \forall(g,g')\in\Ical_G\times\Ical_G,\quad (\T_e)_{g,g'}=\left\{
\begin{aligned}
& \Sigma_{r}^g
\quad\text{ if } g=g',\\
& -\Sigma_{s}^{g'\rightarrow g} \quad\text{ if }g\neq g'{.}
\end{aligned}
\right.
\end{equation*}
The diffusion matrix is  denoted by $\D$. At (almost) every point in $\Omega$, it is a diagonal matrix of  $\R^{G\times G}$:
\begin{equation*}
\forall g\in\Ical_G,\quad \D_{g,g}={ D^g}.
\end{equation*}
Let $\Mbb_f$ be the fission matrix of 
$\R^{G\times G}$ defined for all $(g,g')\in \Ical_G\times\Ical_G$ by
\begin{equation*}
    (\Mbb_f)_{g,g'}=\chi^g(\nu\Sigma_f)^{g'}.
\end{equation*}
Mathematically speaking, the coefficients of the measurable fields of matrices $\T_{e}$ and $\D$ are supposed to be such that:
\begin{equation}\label{MG_Pos}
\left\{
\begin{array}{ll}
(0)&
\forall g,g'\in\Ical_G,\ ({D}^g, \Sigma_{r}^g,\Sigma_{s}^{g'\to g})\in \Pcal W^{1,\infty}(\Omega)\times\Pcal W^{1,\infty}(\Omega)\times L^\infty(\Omega)\,,\\
(i)&\exists\,({D})_*,({D})^*>0,\ \forall\,g\in\Ical_G,\ ({D})_*\le {D}^g\le ({D})^*\mbox{ a.e. in }\Omega\,,\\
(ii)&\exists\,(\Sigma_{r})_*,(\Sigma_{r})^*>0,\ \forall\,g\in\Ical_G,\ (\Sigma_{r})_*\le \Sigma_{r}^g\le (\Sigma_{r})^*\mbox{ a.e. in }\Omega\,,\\
(iii)&\exists {\epsilon\in(0,(G-1)^{-1})},\ \forall\,g,g'\in\Ical_G, g'\neq g,\ |\Sigma_{s}^{g\rightarrow g'}|\leq\epsilon \Sigma_{r}^g\mbox{ a.e. in }\Omega.
\end{array}
\right.
\end{equation}
As a consequence of (\ref{MG_Pos}): the matrix $\T_{e}$ is (almost everywhere) strictly diagonally dominant (so it is invertible)\,; the matrix $\D$ is also invertible (almost everywhere).
\begin{remark}
In general, the removal $\T_e$ is not symmetric due to the scattering cross sections and in the absence of up-scattering, i.e. in the situation where for all $1\leq g<g'\leq G$, we have  $\Sigma_{s}^{g'\rightarrow g}=0$, the matrix $\T_e$ is lower triangular by (energy group) block.
\end{remark}

Let $\phi:=(\phi^1,\ldots,\phi^G)$ and  $\pvec:= (\pvec^1,\ldots,\pvec^G)$ with
$\pvec^g=(p_{x}^g)_{x=1,d}\in\R^d$ for $1\le g\le G$. We also define 
 $\dive\pvec=(\dive_{\xvec} \pvec^g)^{g=1,G}\in\R^G$.
 Then, the multigroup neutron diffusion equation in mixed dual formulation can be written in compact form as:

 \begin{equation}
 \label{eq:mixed-EP}
( ME) \quad
 \left\lbrace
 \begin{aligned}
    \Dbb^{-1}  \pvec + \grad \phi  & =  0 &\text{ in }\Omega,\\
     \dive {\pvec}  + \T_e \phi & = \frac{1}{\keff} \Mbb_f\phi &\text{ in } \Omega, \\
       \phi & = 0 & \text{ on }~\Gamma_D,\\
       \pvec\cdot\nvec & = 0 & \text{ on }~\Gamma_N, \\
       -\,\pvec\cdot\nvec + \frac{1}{2} \phi & = 0 & \text{ on }~\Gamma_R.
 \end{aligned}
 \right.
 \end{equation}

The coefficients of the fission matrix are assumed to be such that
\begin{align}
\left\{
\begin{aligned}
&\forall g \in \Ical_g, (({\nu\Sigma_f})^g,\chi^g)\in L^\infty(\Omega)\times L^\infty(\Omega),\\
&\forall g \in \Ical_g,\, 0\leq ({\nu\Sigma_f})^g \text{ a.e. in }\Omega, \text{ and } (\underline{\nu\Sigma_f})\neq { 0} \text{ a.e. in }\Omega.
\end{aligned}
\right.
 \label{eq:fission_assumption}
\end{align}

Under the assumptions for the coefficients \eqref{MG_Pos} and \eqref{eq:fission_assumption}, the eigenvalue $\keff$ (the multiplication factor) with the greatest modulus is simple, real and strictly positive . Furthermore, the associated eigenfunction $\phi$ (neutron flux) is real and positive, up to a multiplicative constant. We also note that
solving the mixed problem \eqref{eq:mixed-EP} is equivalent to solving the primal problem \eqref{eq:primal-EP}.

\section{The MF-PINNs for the source problem}
\label{section-pinns-SP}
\subsection{Notations}
The shorthand notation $A\lesssim B$ is used for the inequality $A\leq C B$ where $A$ and $B$ are scalar quantities, and $C$ is a generic constant.

Given an open set $\mathcal{O}\subset\R^{{d}}$, ${d}=1,2,3$,
we denote by $\Lvec^2({\mathcal{O}})=(L^2({\mathcal O}))^{{d}}$  and $\Hvec^s(\mathcal{O})= (H^s(\mathcal{O}))^d$ for $s\in \mathbb{R}$.  We use the notation $(\cdot,\cdot)_{0,\mathcal{O}}$ for $L^2({\mathcal O})$ and $\Lvec^2({\mathcal{O}})$ scalar products.
Let $G\in\N \setminus\{0,1\}$ and
given a function space $W$, the product space $W^G$ is denoted by $\Wud$, for instance $\underline{L^2}({\mathcal{O}})$ and $\underline{\Lvec^2}({\mathcal{O}})$. We also extend the notation  $(\cdot,\cdot)_{0,\mathcal{O}}$ for $\underline{L}^2({\mathcal O})$ and $\underline{\Lvec}^2({\mathcal{O}})$ inner products.
 Finally, it
is assumed that the reader is familiar with vector-valued function spaces related to the diffusion equation such as $\Hvec(\dive,\Omega)$. To be convenient, let us define the following function spaces:
\begin{equation*}
    H^1_{0,\Gamma_D}(\Omega):= \left\{ \psi \in H^1(\Omega) \lvert \quad \psi_{|\Gamma_D}=0 \right\},
\end{equation*}
\begin{equation*}
    \Hvec_{0,\Gamma_N}(\dive,\Omega):= \left\{ \qvec \in \Hvec(\dive,\Omega) \lvert \quad {\qvec\cdot\nvec}_{|\Gamma_N}=0 \right\},
\end{equation*}
and 
\begin{equation*}
    \Qvec:=  \left\{ \qvec \in \Hvec_{0,\Gamma_N}(\dive,\Omega) \lvert \quad {\qvec\cdot\nvec}_{|\Gamma_R}\in L^2(\Gamma_R) \right\}.
\end{equation*}

\subsection{The standard loss function for MF-PINNs framework}
First of all, MF-PINNs framework is presented for the source problem. Let us consider a source term $S_f\in \underline{L^2}(\Omega)$, the source problem in mixed dual form of problem \eqref{eq:mixed-EP} can be written as
 \begin{equation}
 \label{eq:mixed-SP}
 \left\lbrace
 \begin{aligned}
    \Dbb^{-1}  \pvec + \grad \phi  & =  0 &\text{ in }\Omega,\\
     \dive {\pvec}  + \T_e \phi & = S_f &\text{ in } \Omega, \\
        \phi & = 0 & \text{ on }~\Gamma_D,\\
       \pvec\cdot\nvec & = 0 & \text{ on }~\Gamma_N, \\
       -\,\pvec\cdot\nvec + \frac{1}{2} \phi & = 0 & \text{ on }~\Gamma_R,
 \end{aligned}
 \right.
 \quad
 \Longleftrightarrow
  \left \lbrace
 \begin{aligned}
   \Lcal [\zeta]  &= S \quad \text{ in }\Omega, \\
  \Bcal[\zeta]  &=0 \quad \text{ on } \partial\Omega,
 \end{aligned}
 \right.
 \end{equation}
 where 
 \begin{equation*}
 \zeta:=(\pvec,\phi), \quad
    \Lcal[\zeta] :=
\begin{pmatrix}
\Dbb^{-1}\pvec   + \grad \phi \\
\dive \pvec + \T_e \phi
\end{pmatrix},
\quad 
 \quad S = 
 \begin{pmatrix}
     0\\
     S_f
 \end{pmatrix},
 \end{equation*}
 and the notation 
 $\Bcal[\pvec,\phi] $ stands for the boundary condition (BC) of the mixed dual form  defined by:
 \begin{equation}
 \label{eq:mixed-BC}
\Bcal[\pvec,\phi] = \left\lbrace
\begin{array}{lll}
     \phi & \text{on} \quad \Gamma_D &\text{ Dirichlet BC}, \\
     \, \pvec\cdot\nvec & \text{on} \quad\Gamma_N  &\text{ Neumann BC}, \\
     -\,\pvec\cdot\nvec + \frac{1}{2} \phi & \text{on} \quad\Gamma_R &\text{ Robin BC}.
\end{array}
\right.
\end{equation}
  Unlike traditional neuron networks with the loss function only coming from data, classical PINNs are unsupervised learning methods that incorporate physical laws directly into the loss function of the neural networks.
  In general, the objective function of the  PINNs framework  is defined by using  the PDE and boundary residual in the $L^2$ norm. Therefore,
 the standard  loss function of MF-PINNs associated to \eqref{eq:mixed-SP} is defined by
  \begin{equation}
  \label{eq:full-loss-MF-PINNs}
\begin{aligned}
    \Jcal[\zeta] 
    & = \lVert \Lcal [\zeta] - S \rVert_{\underline{\Lvec^2}(\Omega)}^2  +\gamma_b \lVert \Bcal[\zeta] \rVert_{\underline{L^2}(\partial \Omega)}^2 \\[0.3cm]
    & = \lVert 	\dive \pvec+ \T_{e}\,\phi - S\rVert_{\underline{L^2}(\Omega)}^2 +  \lVert 	\Dbb^{-1}  \pvec + \grad \phi \rVert_{\underline{\Lvec^2}(\Omega)}^2  +\gamma_b \lVert \Bcal[\zeta] \rVert_{\underline{L^2}(\partial \Omega)}^2, 
\end{aligned}
\end{equation}
where $\gamma_b>0$ is a penalty coefficient.\\
This standard MF-PINNs is studied in \cite{do2025physics} for one group neutron diffusion equation. Then, this strategy is extended to the more general multigroup setting in \cite{do2026MF}.
\begin{remark}
    For heterogeneous coefficients, the loss function of the primal form \eqref{eq:primal-EP} is not well defined since in general $\dive\left(\Dbb  \grad\phi\right) \notin L^2(\Omega)$.
\end{remark}
In fact, the mixed dual form \eqref{eq:mixed-SP} is a first-order system, so the MF-PINNs framework is similar to First Order System Least Squares  (FOSLS). In particular,
let us introduce a bilinear form $b$ and the linear functional F by 
\begin{equation*}
    \begin{aligned}
        b(\zeta,\xi) &:= \displaystyle\left(\Lcal[\zeta],\Lcal[\xi] \right)_{0,\Omega} + \left(\Bcal[\zeta],\Bcal[\xi] \right)_{0,\Gamma_R} \\[0.2cm]
        &=\left( \Dbb^{-1}  \pvec + \grad \phi, \Dbb^{-1}  \qvec + \grad \psi \right)_{0,\Omega} 
        + \left(\dive {\pvec}  + \T_e \phi, \dive {\qvec}  + \T_e \psi \right)_{0,\Omega} \\
         &+ ( -\pvec\cdot\nvec + \frac{1}{2} \phi, -\qvec\cdot\nvec + \frac{1}{2} \psi)_{0,\Gamma_R},
    \end{aligned}
\end{equation*}
and
\begin{equation*}
    F(\xi) := \left(S, \Lcal[\xi]\right)_{0,\Omega}=\displaystyle\left(S_f, \dive {\qvec}  + \T_e \psi \right)_{0,\Omega}.
\end{equation*}
We now consider the following variational formulation: \\
Find $\zeta=(\pvec,\phi)\in \Qvec\times  H^1_{0,\Gamma_D}(\Omega)$ such that
\begin{equation}
\label{eq:variational-problem}
   b(\zeta,\xi) = F(\xi), \quad \forall\xi=(\qvec,\psi)  \in \Qvec\times  H^1_{0,\Gamma_D}(\Omega).
\end{equation}
Following \cite{bernkopf2023optimal,bernkopf2024optimal}, for all $\zeta=(\pvec,\phi)\in \Qvec\times  H^1_{0,\Gamma_D}(\Omega)$,  we have the norm equivalence
\begin{equation}
    \lVert \phi \rVert_{\underline{H}^1(\Omega)}^2+ \lVert \pvec \rVert_{\underline{\Hvec}(\dive,\Omega)}^2 + \lVert \pvec\cdot \nvec \rVert_{\underline{L^2}(\Gamma_R)}^2 \lesssim b(\zeta,\zeta) \lesssim  \lVert \phi \rVert_{\underline{H}^1(\Omega)}^2+ \lVert \pvec \rVert_{\underline{\Hvec}(\dive,\Omega)}^2 + \lVert \pvec\cdot \nvec \rVert_{\underline{L^2}(\Gamma_R)}^2.
\end{equation}
According to Lax-Milgram theorem, the variational formulation \eqref{eq:variational-problem} has unique solution and it is unique minimizer of 
\begin{equation*}
    J(\zeta) = \frac{1}{2} b(\zeta,\zeta) -F(\zeta) = \frac{1}{2}\lVert \Lcal [\zeta] - S \rVert_{\underline{\Lvec^2}(\Omega)}^2 + \frac{1}{2}\lVert -\pvec\cdot\nvec + \frac{1}{2} \phi\rVert_{\underline{L}^2(\Gamma_R)} -\frac{1}{2} \lVert S\rVert_{\underline{\Lvec^2}(\Omega)}^2.
\end{equation*}
Therefore, the unique solution of the mixed dual form \eqref{eq:mixed-SP} satisfies 
\begin{equation*}
    \zeta^{\star} = \displaystyle{\arg\min}_{\zeta \in \Qvec\times  H^1_{0,\Gamma_D}(\Omega)} \left\{ \lVert \Lcal [\zeta] - S \rVert_{\underline{\Lvec^2}(\Omega)}^2 + \lVert -\pvec\cdot\nvec + \frac{1}{2} \phi\rVert_{\underline{L}^2(\Gamma_R)} \right\}.
\end{equation*}
This is equivalent to minimizing the $L^2$ norm of the PDE and boundary residual \eqref{eq:full-loss-MF-PINNs}.
\subsection{Neural Network}
In this work, the aim is to use a neural network as a parametric function with parameter $\theta$  that maps an input $x\in \mathbb{R}^d$ to an output $\zeta_{\theta}(x):=(\pvec_{\theta},\phi_{\theta}) \in \mathbb{R}^{ G(d+1)} $ to approximate the neutron flux and current of the multigroup neutron diffusion equation. In particular,  fully connected neural network (FCNN), the most commonly architecture in PINNs, is considered in this work due to their simplicity and effectiveness in approximating the complex solution. FCNN typically comprises an input layer, $L$ hidden layers, and an output layer as follows
\begin{equation}
\label{eq:NN}
   \begin{aligned}
  h_{\ell} &= \sigma_{\ell}(W_{\ell}h_{\ell-1}+b_{\ell}),  \quad \text{for} \quad h_0 =x \quad \text{and}\quad \ell=1,2,\cdots,L, \\
    \zeta_{\theta}(x)  & = W_{L+1} h_L + b_{L+1},
\end{aligned} 
\end{equation}
  where $\sigma_{\ell}$ stands for an entry-wise activation function of hidden layer $\ell$ such as Tanh and Sin function. The learnable parameters $W_l \in \R^{n_{\ell}\times n_{\ell -1}}$ and $b_l \in \R^{n_{\ell}}$ are respectively weight and bias associated to each layer $\ell$ with the number of neurons $n_{\ell}$. Additionally, let $\theta =\left( W_1,b_1, \cdots, W_L,b_L,W_{L+1},b_{L+1}\right)\in \Theta$ denote the set of weights and biases that are optimized  during the training process. 
\subsection{Hard boundary condition for MF-PINNs}
 Classically,  the boundary condition term is incorporated into the Pinns loss function \eqref{eq:full-loss-MF-PINNs}  as a penalty term  and this method is commonly referred as the \textit{soft boundary condition}.
However, a neural network using \textit{soft boundary condition} method does not exactly satisfy the boundary conditions. Training perfectly the boundary term is  very challenging and  it may lead to the imbalance gradients problem of the training process. In order to guaranty the boundary condition,  one can build a neural network in such a way that $\Bcal[\zeta_{\theta}] = 0 \quad \text{on }  \partial \Omega$ which helps to avoid the imbalance gradient problem, as well as to improve the convergence of the optimization problem. In fact, one can introduce an extra modified layer as:
\begin{equation*}
\tilde{\zeta}_{\theta}(x) := F_b(x,\zeta_{\theta}(x)),
\end{equation*}
where $F_b: \R^{d}\times  \R^{G(d+1)} \rightarrow  \R^{G(d+1)} $ is a $C^1$ function that depends on the boundary function $b(x)$ and $\Omega$, which is constructed
to modify the output function for the boundary condition.  This method is commonly referred as hard boundary conditions (HBC). 
Hard enforcement Dirichlet boundary condition can be done by modifying the output of neural network by 
\begin{equation*}
    \tilde{\phi}_{\theta}(x)= L(x)\phi_{\theta}(x) +B(x),
\end{equation*}
where $L(x)$ is a distance function to the Dirichlet boundary and $B(x)$ is an extension of $b(x)$ such that 
$B(x)=b(x)$ at the Dirichlet boundary. It is also interesting to point out that
for square or rectangular domains ($\nvec =(\pm 1,0)$ or $\nvec =(0,\pm 1)$), which are common in nuclear reactors, reflective boundary condition ($\pvec\cdot \nvec =0 $)
or vacuum boundary condition ($\pvec\cdot \nvec = \frac{1}{2}\phi$) becomes Dirichlet boundary condition on the appropriate current component. For arbitrary geometries and the general framework for enforcing hard boundary conditions in PINNs, the reader is referred to \cite{xie2023automatic,xie2024boundary,lyu2021enforcing, bersetche2023deep}.

\subsection{The Physic-Based Loss Scaling method}
\begin{figure}[H]
      \centering
      \includegraphics[scale=0.95]{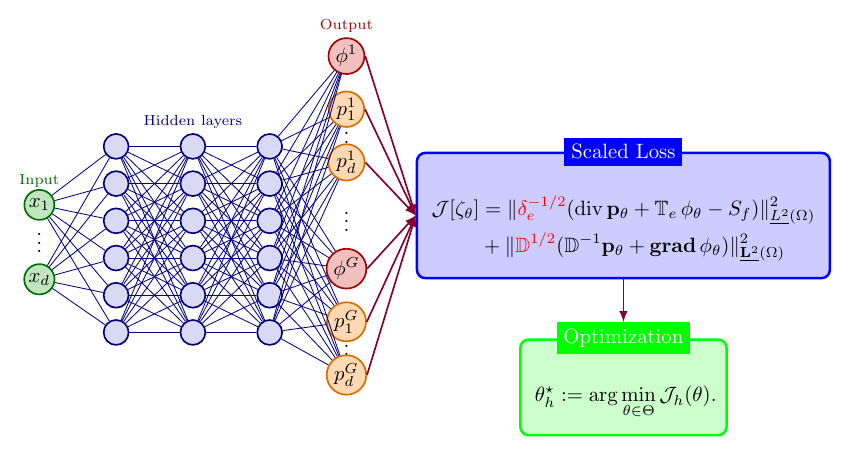}
      \caption{ Fully connected neural network  with scaled loss  for MF-PINNs}
      \label{fig-MF-PINNs}
  \end{figure}

In~ \cite{ciarletChapter5}, some of the co-authors propose a posteriori \textit{residual} estimates for the mixed finite element discretization of the multigroup neutron diffusion equation. These estimates involves a so-called reconstruction denoted $\tilde{\zeta}_h=(\tilde{\pvec}_h,\tilde{\phi}_h)\in\Qvec\times H^1_{0,\Gamma_D}(\Omega)$. Let $\Tcal_h$ be a mesh and $K\in\Tcal_h$ be a mesh element. In the case where $\Gamma_R=\emptyset$, the local a posteriori estimator writes
\begin{align}
   \eta_K:= \left({\eta}^2_{r,K} +\sum_{K'\in N(K)}\eta^2_{f,K'}\right)^{1/2},
\end{align}
where \begin{align*}
N(K)&=\{K'\in\Tcal_h\ |\ \text{dim}_{H}({\partial K'\cap \partial K})=d-1\}, \\
    {\eta_{r,K}}&=\|\delta_e^{-1/2} (\dive \tilde{\pvec}_h +\T_e\tilde{\phi}_h-S_f)\|_{\underline{\Lvec^2}(K)}\, ,\\
    {\eta_{f,K}} &= \|\Dbb^{1/2}(\Dbb^{-1}\tilde{\pvec}_h+\grad \tilde{\phi}_h)\|_{\underline{\Lvec^2}(K)}\, ,
\end{align*}
with dim$_H$ the Hausdorff dimension, $\delta_{e}$ the diagonal part of $\T_{e}$, and $\delta_{e}^{1/2}$ its square root. In~\cite{CDGM26}, some of the co-authors analyse the case where $\Gamma_R\neq\emptyset$, and a specific estimator is designed to handle the Robin boundary condition. 
Let us note that the hard boundary condition (HBC) is applied  for MF-PINNs. Hence, the boundary conditions are enforced directly in the neural network construction, and the boundary loss term is not included in the total loss function.
Then, in order to accelerate the training of MF-PINNs, we propose the following physics based loss scaling for MF-PINNs:
 \begin{equation}
 \label{eq:PBLS}
    \Jcal[\zeta] 
   := \lVert \delta_e^{-1/2}	(\dive \pvec+ \T_{e}\,\phi - S_f)\rVert_{\underline{L^2}(\Omega)}^2 +  \lVert 	\Dbb^{1/2}(\Dbb^{-1}  \pvec + \grad \phi) \rVert_{\underline{\Lvec^2}(\Omega)}^2.
\end{equation}
  The MF-PINNs framework with scaled loss function is presented in Figure \ref{fig-MF-PINNs}.

\begin{remark}
Under the assumptions for the coefficients \eqref{MG_Pos},
     the scaled loss function defined in \eqref{eq:PBLS} is actually an equivalent norm to the unscaled loss function (classical $L^2$ norm).
\end{remark}
\begin{remark}
     The scaled loss function  \eqref{eq:PBLS} can be seen as
    \begin{equation*}
    \Jcal[\zeta] 
   := \lVert \mathbb{M}^{1/2}	(\dive \pvec+ \T_{e}\,\phi - S_f)\rVert_{\underline{L^2}(\Omega)}^2 +  \lVert 	\mathbb{N}^{1/2}(\Dbb^{-1}  \pvec + \grad \phi) \rVert_{\underline{\Lvec^2}(\Omega)}^2, 
\end{equation*}
   where $\mathbb{M}= \delta_e^{-1}$ and $\mathbb{N}=\Dbb$ are two preconditioning operators which are positive define and symmetric.
\end{remark}

\subsection{Training MF-PINNs}

For any approximate solution of a neural network  $\zeta_{\theta}=(\pvec_{\theta},\phi_{\theta})$.
  Based on the loss scaling method discussed above, the objective function for MF-PINNs is given by
 \begin{equation}
 \label{eq:scaling-loss-mixed-SP}
\begin{aligned}
    \Jcal[\zeta_{\theta}] 
   := \lVert \delta_e^{-1/2}	(\dive \pvec_{\theta}+ \T_{e}\,\phi_{\theta} - S_f)\rVert_{\underline{L^2}(\Omega)}^2 +  \lVert 	\Dbb^{1/2}(\Dbb^{-1}  \pvec_{\theta} + \grad \phi_{\theta}) \rVert_{\underline{\Lvec^2}(\Omega)}^2.  
\end{aligned}
\end{equation}
 Instead of working with the continuous loss function \eqref{eq:scaling-loss-mixed-SP}, it is common to approximate it by the following empirical loss:
\begin{equation}
\label{eq:emp-mixed-SP}
\begin{aligned}
    \Jcal_{h}(\zeta_{\theta}) &:= \ds \frac{1}{N_{r}}\sum_{i=1}^{N_r} \left[ \Lcal[\zeta_{\theta}](x_i^r) - S(x_i^r)\right]^2 \\
    &= \ds \frac{1}{N_{r}}\sum_{i=1}^{N_r} \delta_e^{-1}(x_i^r) (\dive \pvec_{\theta}+ \T_{e}\,\phi_{\theta} - S_f)^2(x_i^r) \\
    &+ \ds \frac{1}{N_{r}}\sum_{i=1}^{N_r}\D(x_i^r) (\Dbb^{-1}  \pvec_{\theta} + \grad \phi_{\theta})^2(x_i^r),
\end{aligned}
\end{equation}
where $N_r$ stands for the number of residual points,
$\left\{ x_i^r\right\}_{i=1}^{N_r}$ are \textit{residual points} randomly sampled in $\Omega$. It is important to note that the sampling method (SM) for the residual points in the empirical loss \eqref{eq:emp-mixed-SP} plays an important role in the performance of PINNs. \\
Finally, we have  to minimize the discrete loss function
\begin{equation*}
    \theta^{\star}_h := \arg\min_{\theta\in \Theta} \Jcal_h(\theta).
\end{equation*}
This problem can be solved with gradient-based optimization methods such as the gradient descent (GD):
\begin{equation*}
 \theta_{j+1} = \theta_j -\eta_j  \nabla \Jcal_h(\theta_j),    \quad j=1,\ldots,J,
\end{equation*}
where $J$ is the number of inner iterations, $\theta_j$ represents the parameters of the neural network at iteration $j$ and $\eta_j$ is the learning rate. However, the gradient descent method  often suffers from slow convergence and sensitivity to the choice of learning rate. To address these limitations, Adam (Adaptive Moment Estimation) optimizer \cite{Kingma2014AdamAM} combines the advantages of the momentum and adaptive learning rate method to  improve convergence speed and stability. Therefore, Adam optimizer is widely used in the  training of neural networks, including PINNs.

\section{The MF-PINNs for the k-eigenvalue problem}
\label{section-MF-PINNs-EP}
In this section, we would like to use the MF-PINNs framework with scaled loss function to approximate the solution of the k-eigenvalue problem. In fact, it is still a challenge  to accurately predict the eigenvalue and the associated eigenfunction of an eigenvalue problem using PINNs. Currently,
there are two main  approaches for PINNs applied to the k-eigenvalue problem. The first  uses the multiplication factor $\keff$ as a trainable parameter, similar to the parameters of a neural network. This value is  optimized simultaneously with the parameter $\theta$ of the neural network during the training process. By this way, the eigenvalue $\keff$ is learned by minimizing the physical loss function. 
Although this method is simple to implement in the PINNs framework,  it may converge to trivial or inaccurate solutions. Therefore, this approach usually requires extra physical constraints to avoid trivial solutions such as the regularization method in \cite{elhareef2023physics}, fixed point constraint PINNs (FC-PINNs) \cite{bi2025fc} and amplitude constraint \cite{yoo2025physics}. In the second approach,  the value of $\keff$ is calculated using an iterative procedure based on classical fixed point methods such as power iteration (PI).
In the symmetric case, the neural network only learns to predict the eigenfunction and the eigenvalue is updated externally using the Rayleigh quotient method  in~\cite{yang2023neural,yang2023physics}. In general, this  approach is more stable since the calculation of $\keff$  mimics traditional numerical methods but  clearly introduces an additional outer loop. Obviously, it increases the computational complexity of the PINNs. Furthermore, the convergence of this method  is very slow if the dominance ratio is high (close to one), as in several reactor physics applications. Some acceleration methods such as Anderson acceleration \cite{walker2011anderson}  and momentum acceleration \cite{xu2018accelerated} can be applied to reduce the number of outer iterations \cite{calloo2023anderson}. These acceleration methods can be applied to  improve the convergence of MF-PINNs \cite{do2025physics}.

In this work, the k-eigenvalue problem is addressed using the inverse power iteration described below. Starting from an initial guess for the eigenvalue  $\keff^0$ and the eigenfunction $\zeta^0_{\theta}=(\pvec^0_{\theta}, \phi^{0}_{\theta})$,  at each iteration $n+1$, we compute the updated eigenpair $(\keff^{n+1},\zeta^{n+1}_{\theta})$. In particular, given $(\keff^{n},\zeta^{n}_{\theta})$,   we can first obtain $\zeta^{n+1}_{\theta}$ by solving the associated fixed source problem in mixed form:
 \begin{equation}\label{eq:IP-source}
         \left\lbrace
         \begin{aligned}
            \Dbb^{-1}  \pvec^{n+1}_{\theta} + \grad \phi^{n+1}_{\theta}  & =  0 &\text{ in }\Omega,\\
             \dive {\pvec^{n+1}_{\theta}}  + \T_e \phi^{n+1}_{\theta} & = \frac{1}{\keffn}S_f^n &\text{ in } \Omega, \\
               \Bcal [\pvec^{n+1}_{\theta},\phi^{n+1}_{\theta}] & = 0 & \text{ on } \partial\Omega,
         \end{aligned}
         \right.
         \end{equation}
where   $S_f^{n}=\Mbb_f \phi_{\theta}^{n} $.      
In fact, this source problem is solved by using MF-PINNs with the Physics-Based Loss Scaling method until the fixed number of inner iterations denoted by $J$. Then, the value of $\keff$ and the source term are updated at each inverse power iteration (outer iteration) as 
\begin{equation*}
  S_f^{n+1}=\Mbb_f \phi_{\theta}^{n+1} \quad\text{and} \quad
		\keffnp = \keffn \displaystyle\frac{[S_f^{n+1}(\xvec^r)]^{T}[S_f^{n+1}(\xvec^r)]}{ [S_f^{n+1}(\xvec^r)]^{T}[S_f^{n}(\xvec^r)] }. 
\end{equation*}
This process is repeated until convergence is achieved. We also note that
it is unnecessary to construct the fission matrix $\Mbb_f$ in practice
since 
\begin{equation*}
  \Mbb_f \phi_{\theta}^{n+1}=(\chi^1 f^{n+1}, \ldots, \chi^G f^{n+1}) 
  \text{ where }  f^{n+1} = \displaystyle\sum_{g'=1}^G \nu\Sigma_f^{g'} \phi_{\theta}^{g',n+1}.
\end{equation*}
Therefore, we can also update the value of $\keff$ in the lower dimension space by 
\begin{equation*}
    \keffnp =\keffn \frac{[f^{n+1}(\xvec^r)]^{T}[f^{n+1}(\xvec^r)]}{[f^{n+1}(\xvec^r)]^{T}[f^{n}(\xvec^r)]}.
\end{equation*}

The details of the MF-PINNs framework is applied to the k-eigenvalue problem is described in Algorithm \ref{al:IP}.

\begin{algorithm}[H]
	\caption{The Inverse Power  Iteration for MF-PINNs with scaled loss function}
    \label{al:IP}
\KwIn{ The  threshold numbers
		 $ \Ecal_{\phi}$ and $ \Ecal_{\keff}$}
\KwOut{ The solution $\zeta^n_{\theta}=(\pvec^n_{\theta}, \phi^{n}_{\theta})$ and $ \keffn$ at iteration number $n$.   }
Set initial guess $\keff^{0}$ and $S_f^0$, $n=0$\\
		\Do{ $\epsilon_{\phi}^{n}\leq \Ecal_{\phi}$ \text{ and }  $ \epsilon_{\keff}^{n} \leq {\Ecal_{\keff }  } $ }{
    \For{$j=1,2,\cdots,J$}{
    
    Calculate discrete scaled loss in mixed form:
    \begin{equation*}
        \begin{aligned}
        \Jcal_h(\zeta_{\theta_j^n})
       &= \ds \sum_{i=1}^{N_r} \delta_e^{-1}(x_i^r)	\left(\dive \pvec_{\theta_j^n}+ \T_{e}\,\phi_{\theta_j^n} - S_f^n \right)^2(x_i^r) \\
       &+ \ds \sum_{i=1}^{N_r}  	\Dbb(x_i^r)\left(\Dbb^{-1} \pvec_{\theta_j^n} + \grad \phi_{\theta_j^n}\right)^2(x_i^r).   
    \end{aligned}
    \end{equation*}
    
    Update parameter $\theta$ via gradient-based method:
    \begin{equation*}
	\ds \theta_{j+1}^n = \theta_j^n -\alpha_j \nabla_{\theta}\Jcal_h(\zeta_{\theta_j^n});
	\end{equation*}
	}
		Update source and $\keff$:   
        \begin{equation*}
           S_f^{n+1}=(\chi^1 f^{n+1}, \ldots, \chi^G f^{n+1}) \text{ and }  \ds \keffnp  =\keffn \frac{[f^{n+1}(\xvec^r)]^{T}[f^{n+1}(\xvec^r)]}{[f^{n+1}(\xvec^r)]^{T}[f^{n}(\xvec^r)]},
        \end{equation*}
        where $f^{n+1} = \displaystyle\sum_{g'=1}^G \nu\Sigma_f^{g'} \phi_{\theta_{J}^n}^{g'}$\;
	    \indent
		Evaluate residuals:
        \begin{equation*}
          \epsilon_{\phi}^{n+1}= \displaystyle\frac{\|\phi^{n+1}_{\theta}(\xvec^r)-\phi^{n}_{\theta}(\xvec^r)\|_2}{\|\phi^{n}_{\theta}(\xvec^r)\|_2}
	 \text{ and } \epsilon_{\keff}^{n+1}= \displaystyle\frac{|\keffnp-\keffn|}{\keffn}.  
        \end{equation*}
		
		$n\leftarrow n+1$
		
	}
\end{algorithm}

\section{Numerical results}\label{section:numerical}
For quantitative comparison purposes, let us define the absolute error of $\keff$ as follows:
\begin{equation*}
\Delta \keff = 10^5\times\lvert \keff -\keffref \rvert ~ (pcm),
\end{equation*}
and the relative error for the neutron flux and neutron current are,  respectively, given by
\begin{equation*}
    \bar{\Delta}{\phi}= 10^2 \times \frac{\lVert \phi_{\theta}(\xvec_{\rm test})-\phi_{\rm ref}(\xvec_{\rm test})\rVert_2}{\lVert \phi_{\rm ref}(\xvec_{\rm test})\rVert_2} ~(\%),
\end{equation*}
and 
\begin{equation*}
 \bar{\Delta}{\pvec}= 10^2 \times \frac{\lVert \pvec_{\theta}(\xvec_{\rm test})-\pvec_{\rm ref}(\xvec_{\rm test})\rVert_2}{\lVert \pvec_{\rm ref}(\xvec_{\rm test})\rVert_2} ~ (\%),
\end{equation*}
where $\lVert\cdot\rVert_2$ is the Euclidean norm or the 2-norm, $\xvec_{\rm test} =\left\{ x_i^t\right\}_{i=1}^{N_t}$ is the set of test points. The calculation of the reference solution is performed using the MINOS solver~\cite{baudron2007minos}
of the deterministic code APOLLO3\textsuperscript{\textregistered} \cite{mosca2024apollo3}.

In the previous work \cite{do2025physics},  the MF-PINNs are investigated using 
different fully connected neural network (FCNN) \eqref{eq:NN} architectures with varying depths and widths ($L=3,5,7$ and $N_{\ell}=32,64$). In general, deeper neural networks can improve the accuracy of MF-PINNs, but they also increase the computational cost and the risk of overfitting. Therefore, in this work, we stay with an FCNN  using $L=5$ hidden layers and $N_{\ell}=64$ neurons for each layer as a compromise between the accuracy and computational efficiency.  The hard boundary condition (HBC) is imposed to the neural network, so there is no boundary term in the loss function. Adam optimizer and multi-step learning rate  is applied to the training process of both source and eigenvalue problem. For the source problem, the initial learning rate is   $\eta_0 = 10^{-3}$ and the  decay rate is $\gamma=0.05$ for each $2000$ training iteration.
For k-eigenvalue problem, the initial learning rate is   $\eta_0 = 2\times 10^{-4}$ and the decay rate is $\gamma=0.1$ for each $10000$ training iteration. Moreover, the number of inner iterations is fixed at $J=2000$ and the stopping
criteria are set to $ \mathcal{E}_{\phi} = 10^{-5}$ and $ \mathcal{E}_{\keff}=10^{-6}$. 
We show only the results for the random sampling method with Tanh activation (classical sampling and activation) and the Sobol sampling with Sin activation, which is the best option for MF-PINNs  from our previous studies \cite{do2025physics,do2026MF}.

\subsection{Test cases for the fixed source problem}
\subsubsection{The IAEA pool type reactor benchmark.}
\begin{figure}[H]
	\centering
	\subfloat[ ]{\includegraphics[scale=0.6]{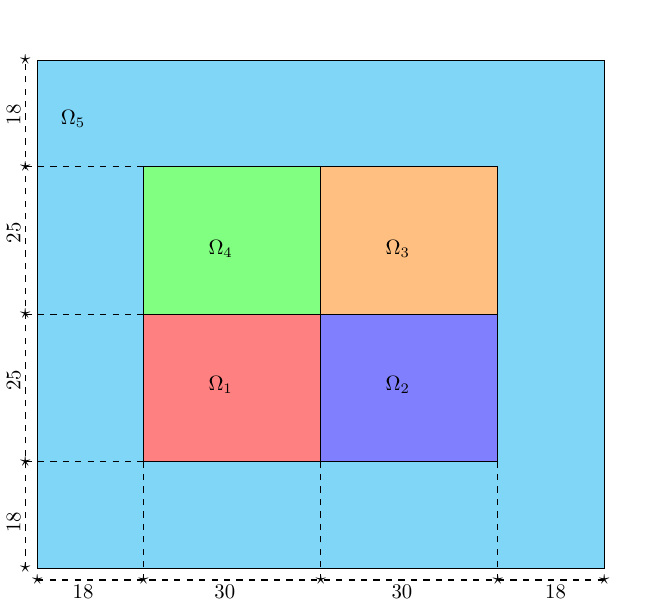}}	
	\caption{ The geometry of the IAEA pool  reactor problem.}
	\label{fig:geo-5-region} 
\end{figure}

We now investigate the IAEA pool type reactor benchmark for one energy group studied in \cite{stepanek1982calculation,cao2020deterministic}  with a total 5 regions corresponding to 5 materials. The geometry of this benchmark depicted in  Figure \ref{fig:geo-5-region}  is  a rectangle of size $96\text{cm}\times 86\text{cm}$, consisting of  two large fuel zones ($\Omega_1$ and $\Omega_3$) and two large absorber zones ($\Omega_2$ and $\Omega_4$) of
size $30\text{cm} \times 25\text{cm}$  for each zone,
surrounded by light water ($\Omega_5$). The  Vacuum boundary condition is applied in this test and the material information is listed in Table \ref{tab:data-CR}. The  reference solution is calculated on the uniform mesh with a cell size of $1.0\text{cm}\times 1.0 \text{cm}$ and the number of residual points is $N_r=10240$.

\begin{table}[H]
  \centering 
     \caption{Cross sections of the IAEA pool reactor problem} 
  \label{tab:data-CR}
  \begin{tabular}{@{}c ccccc@{}}
  \toprule[\heavyrulewidth]
& $\Omega_1$ & $\Omega_2$ & $\Omega_3$ &  $\Omega_4$ & $\Omega_5$ \\ 
  \midrule
  $\Sigma_t $       & 0.60 &   0.48 &  0.70 &  0.65 & 0.90\\
  $\Sigma_s $        & 0.53 &  0.20 &  0.66 &  0.50 & 0.89\\
  $S_f$              & 1.0 &   0.0 &   1.0 &   0.0  & 0.0\\
  \bottomrule[\heavyrulewidth] 
  \end{tabular}
 
\end{table}

\begin{figure}[H]
	\centering
	\subfloat[$\bar{\Delta} \phi$]{\includegraphics[scale=0.4]{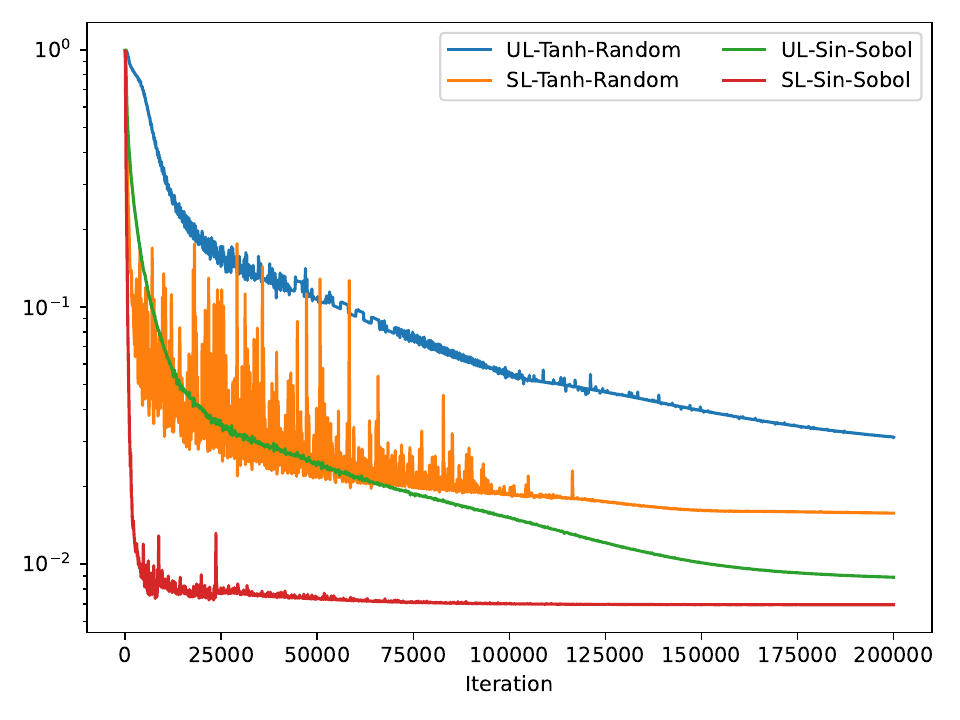}}
    \subfloat[$\bar{\Delta} \pvec$]{\includegraphics[scale=0.4]{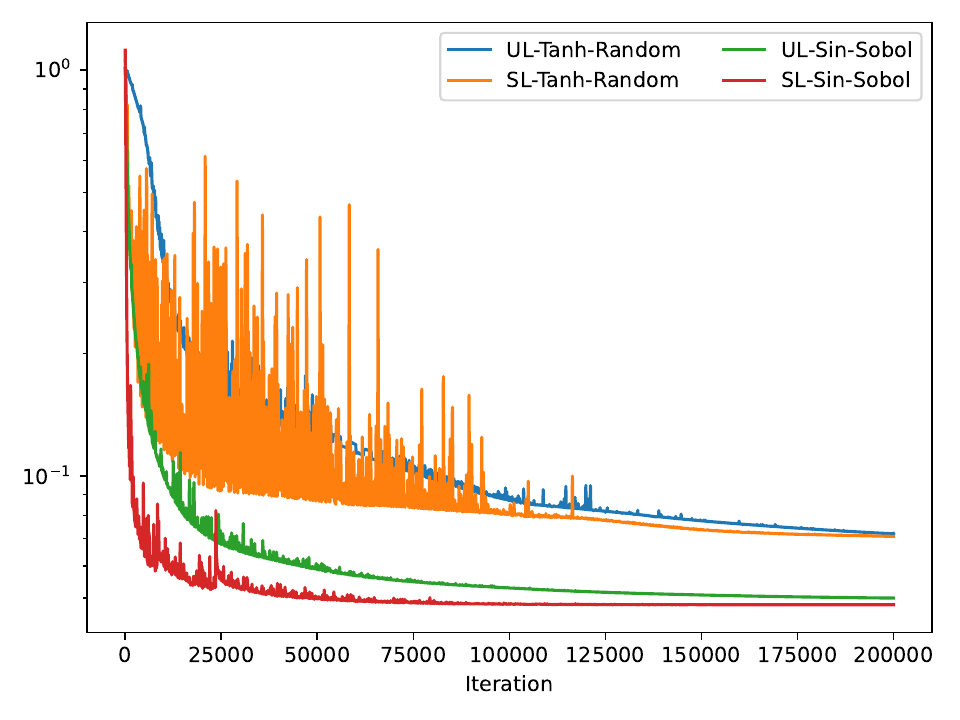}}	
	\caption{ The relative error of the neutron flux and current of the  one-group IAEA pool type reactor}
    \label{fig:IAEA-pool-error}
\end{figure}

\begin{figure}[H]
	\centering
	\subfloat[]{\includegraphics[scale=0.3]{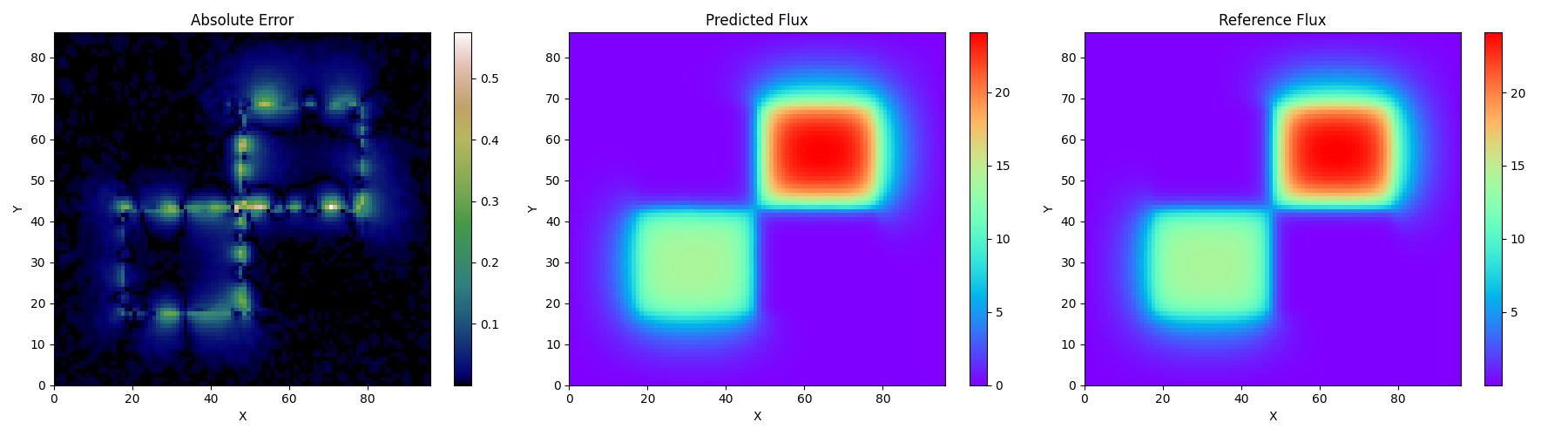}}\\	
	\caption{ Neutron flux and the deviation from the reference solution for the IAEA pool type reactor benchmark with the scaled loss function,  Sobol sampling method and $Sin$ activation function.}
	\label{fig:IAEA-pool-flux} 
\end{figure}

As illustrated in Figure \ref{fig:IAEA-pool-error}, the Sobol sampling method combined with the Sin activation function clearly outperforms the random sampling method with the Tanh activation function. More importantly, the scaled loss (SL) using the PBLS method provides higher accuracy and significantly faster convergence than standard unscaled loss (UL). It is also observed that the prediction error in the scalar flux is lower than that in the neutron current. Figure \ref{fig:IAEA-pool-flux} shows that  MF-PINNs with scaled loss provide a good approximate solution to the IAEA pool type reactor and there are small errors in the interactions between materials.

\subsubsection{The simplified C5G7 benchmark source problem.}
\begin{figure}[H]
\centering
\begin{minipage}{0.5\textwidth}
  \centering
  \includegraphics[scale=0.7]{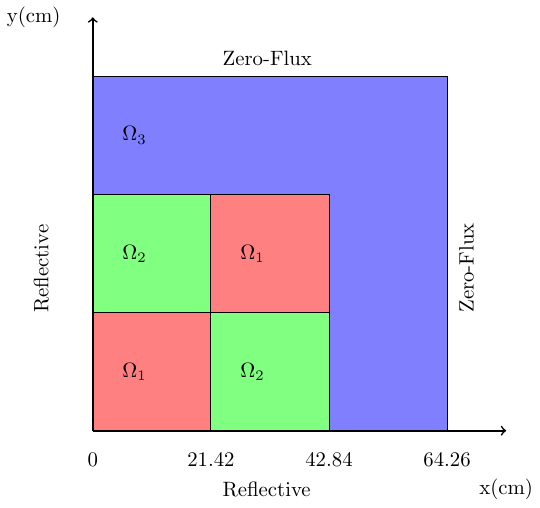} 
  \caption{The simplified C5G7 Geometry}
  \label{fig:C5G7S-geo}
\end{minipage}\hfill
\begin{minipage}{0.5\textwidth}
  \centering
  \captionof{table}{Cross sections $(cm^{-1})$ for the simplified C5G7 benchmark source problem }
  \begin{tabular}{@{}c ccc@{}}
  \toprule
   & $\Omega_1$ & $\Omega_2$ & $\Omega_3$  \\ 
  \midrule
  $D^1 $       & 1.2 &   1.2 &  1.2 \\
  $D^2 $       & 0.4 &   0.4 &  0.2 \\
\addlinespace
  $\Sigma_{r}^1 $       & 0.03 &   0.03 &  0.051 \\
  $\Sigma_{a}^2 $       & 0.3 &   0.25 &  0.04 \\
\addlinespace
  $\Sigma_{s}^{1\to2} $        & 0.015 &  0.015 &  0.05 \\
\addlinespace
  $S_{f}^1 $               & 0.0075 &   0.0075 &   0 \\
  $S_{f}^2 $               & 0 &   0.0 &   0 \\
  \bottomrule
  \end{tabular}
  \label{tab:C5G7S-xs}
\end{minipage}
\end{figure}

 We examine a simplified version of the C5G7 benchmark, originally  studied in \cite{elhareef2023physics} for the k-eigenvalue problem, with homogenized cross sections and two energy groups. In this section, we modify it to a fixed source problem \eqref{eq:mixed-SP}. The geometry of this test case is presented in Figure \ref{fig:C5G7S-geo} with a total of 3 materials where the left and bottom sides are  reflective boundary conditions while the right and top sides are zero flux boundary conditions.
 The cross sections data are provided in Table \ref{tab:C5G7S-xs}. The calculation of the reference solution is performed on a uniform 120 × 120 grid using the MINOS solver and the number of residual points for training is $N_r=20480$.
 
\begin{figure}[H]
	\centering
	\subfloat[$\bar{\Delta} \phi$]{\includegraphics[scale=0.4]{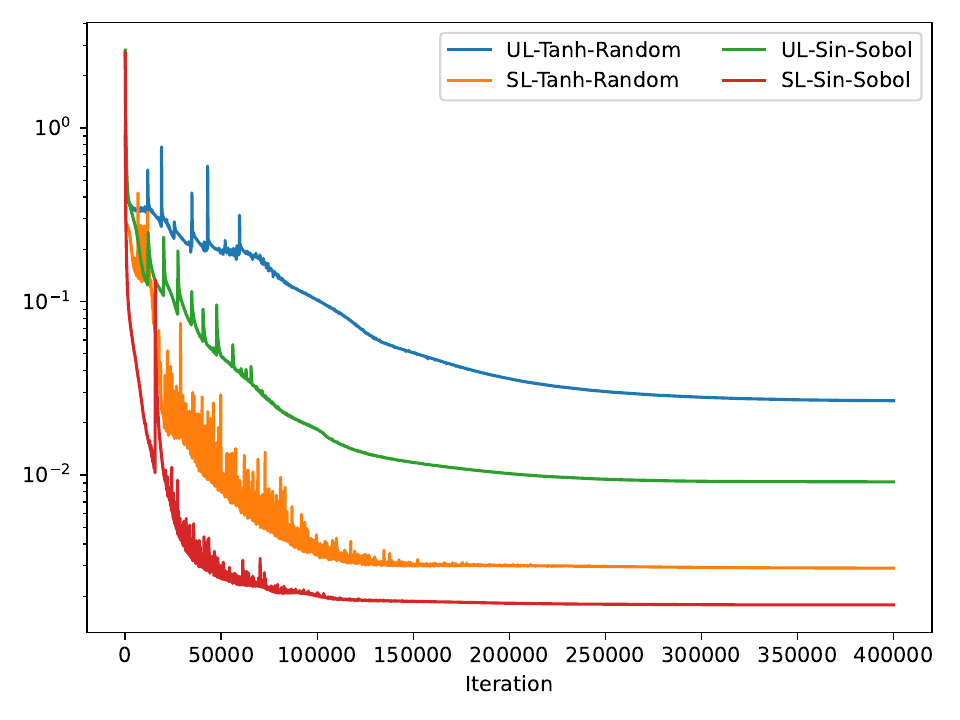}}
    \subfloat[$\bar{\Delta} \pvec$]{\includegraphics[scale=0.4]{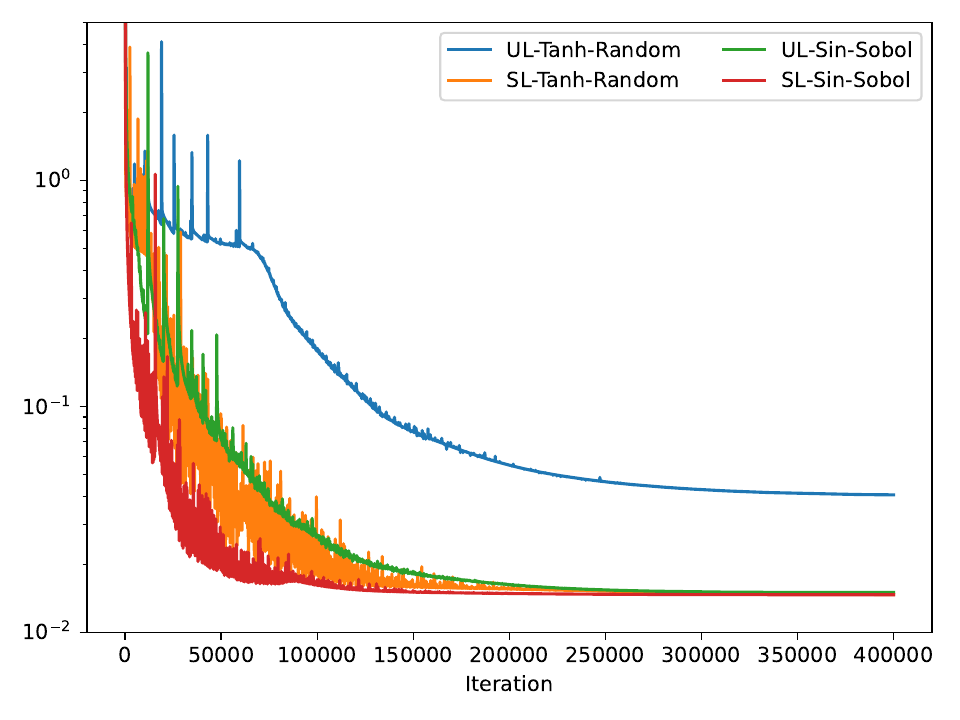}}	
	\caption{ The relative error of the neutron flux and current of the  simplified C5G7 source problem.}
    \label{fig:C5G7S-error}
\end{figure}

\begin{figure}[H]
	\centering
	\subfloat[Fast group ]{\includegraphics[scale=0.3]{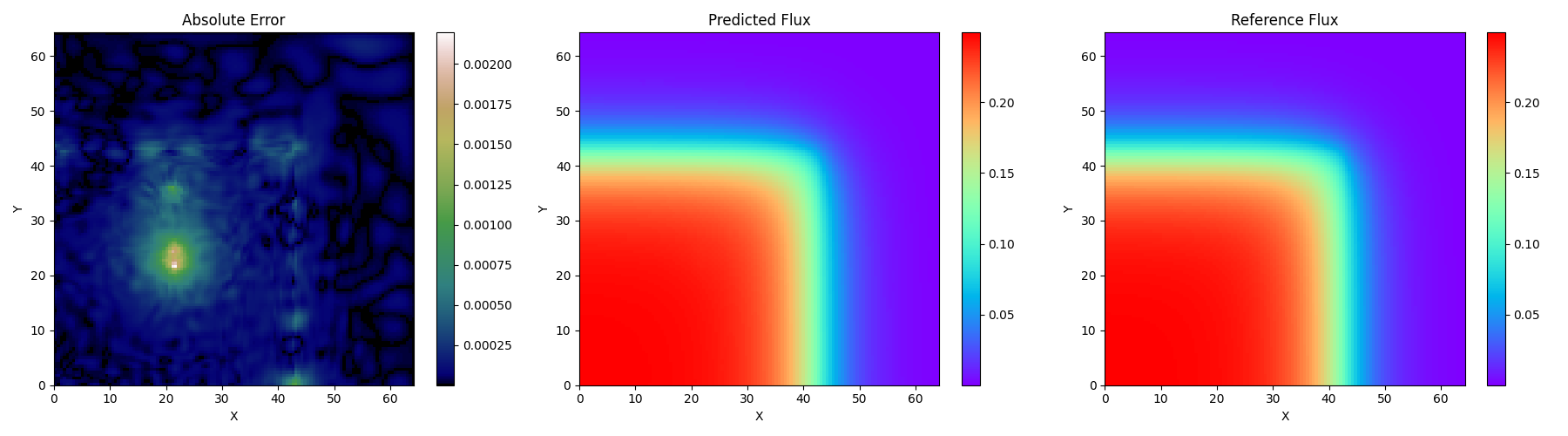}}\\	
	\subfloat[ Thermal group ]{\includegraphics[scale=0.3]{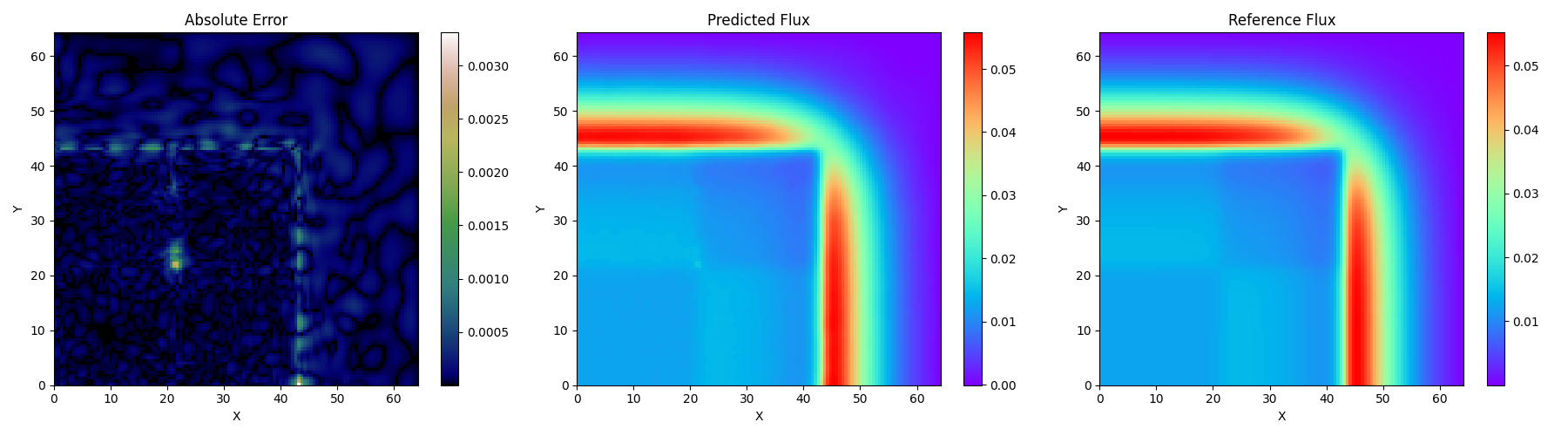}}\\
	
	\caption{ Neutron flux and the deviation from the reference solution for the simplified C5G7 benchmark source problem with scaled loss function, Sobol sampling method and $Sin$ activation function.}
	\label{fig:C5G7S-flux} 
\end{figure}

According to Figures \ref{fig:C5G7S-error} and \ref{fig:C5G7S-flux} the multigroup case exhibits the same behavior as in the single-group case. The Sobol sampling method combined with Sin activation function still provides better performance than the other, and the scaled loss (SL) leads to faster convergence and higher accuracy than the unscaled loss (UL).  MF-PINNs with  scaled loss (SL) accurately approximate both the fast and thermal neutron fluxes. This confirms that  the good results of MF-PINNs with scaled loss for the one group case (IAEA pool type reactor) are still valid for the multigroup case (the simplified C5G7 benchmark).

\subsection{Test cases for the k-eigenvalue problem}
In this numerical section,  we first present all test configurations for the k-eigenvalue problem. The analysis of the obtained results is then discussed at the end of the section since all test cases exhibit similar behavior.
\subsubsection{The simplified C5G7 benchmark}
\begin{figure}[H]
\centering
\begin{minipage}{0.5\textwidth}
  \centering
  \includegraphics[scale=0.7]{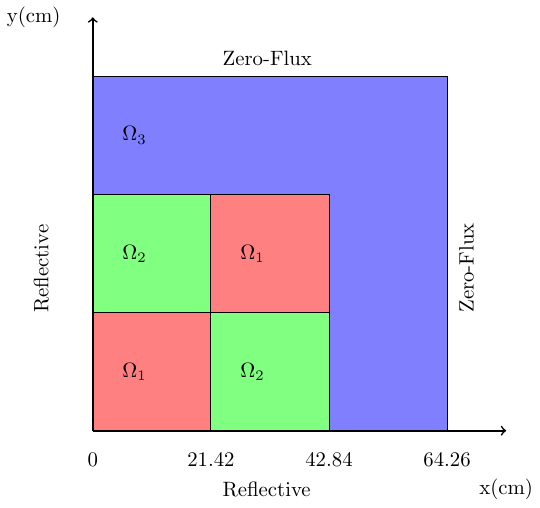} 
  \caption{The simplified C5G7 Geometry}
  \label{fig:C5G7-geo}
\end{minipage}\hfill
\begin{minipage}{0.5\textwidth}
  \centering
  \captionof{table}{Cross sections $(cm^{-1})$ for the simplified C5G7 benchmark }
  \begin{tabular}{@{}c ccc@{}}
  \toprule
 & $\Omega_1$ & $\Omega_2$ & $\Omega_3$  \\ 
  \midrule
  $D^1 $       & 1.2 &   1.2 &  1.2 \\
  $D^2 $       & 0.4 &   0.4 &  0.2 \\
\addlinespace
  $\Sigma_{r}^1 $       & 0.03 &   0.03 &  0.051 \\
  $\Sigma_{a}^2 $       & 0.3 &   0.25 &  0.04 \\
\addlinespace
  $\Sigma_{s}^{1\to 2} $        & 0.015 &  0.015 &  0.05 \\
\addlinespace
  $\nu\Sigma_{f}^1 $               & 0.0075 &   0.0075 &   0 \\
  $\nu\Sigma_{f}^2 $               & 0.45 &   0.375 &   0 \\
  \bottomrule
  \end{tabular}
  \label{tab:C5G7-xs}
\end{minipage}
\end{figure}

 We now consider the simplified version of the C5G7 benchmark (SC5G7) for the  k-eigenvalue problem \cite{elhareef2023physics}. The geometry of this test case is presented in Figure \ref{fig:C5G7-geo} and the cross sections data are provided in Table \ref{tab:C5G7-xs}. The calculation of the reference solution is performed on a uniform 120 × 120 grid using the MINOS solver and the reference value of the multiplication factor is $\keffref= 0.92757$. The number of residual points is $N_r=20480$.

\begin{figure}[H]
	\centering
	\subfloat[$\bar{\Delta} \phi$]{\includegraphics[scale=0.4]{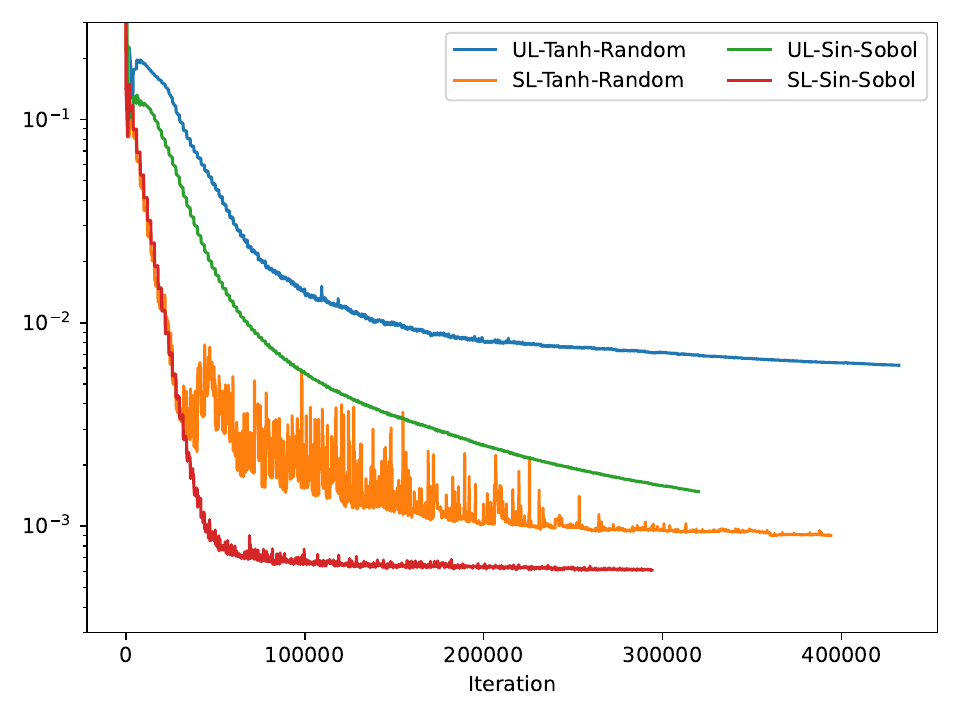}}
    \subfloat[$\bar{\Delta} \pvec$]{\includegraphics[scale=0.4]{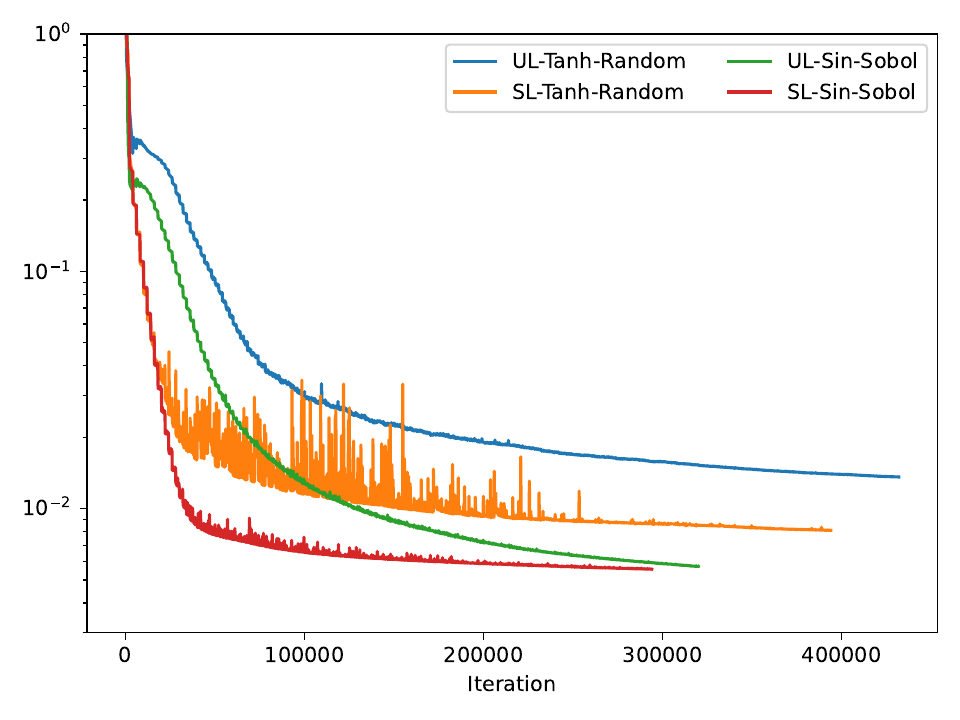}}	\\
	\caption{ The relative error of the neutron flux and  neutron current  of the  simplified C5G7 problem.}
    	\label{fig:C5G7-error} 
\end{figure}

\begin{figure}[H]
	\centering
	\subfloat[Fast group ]{\includegraphics[scale=0.3]{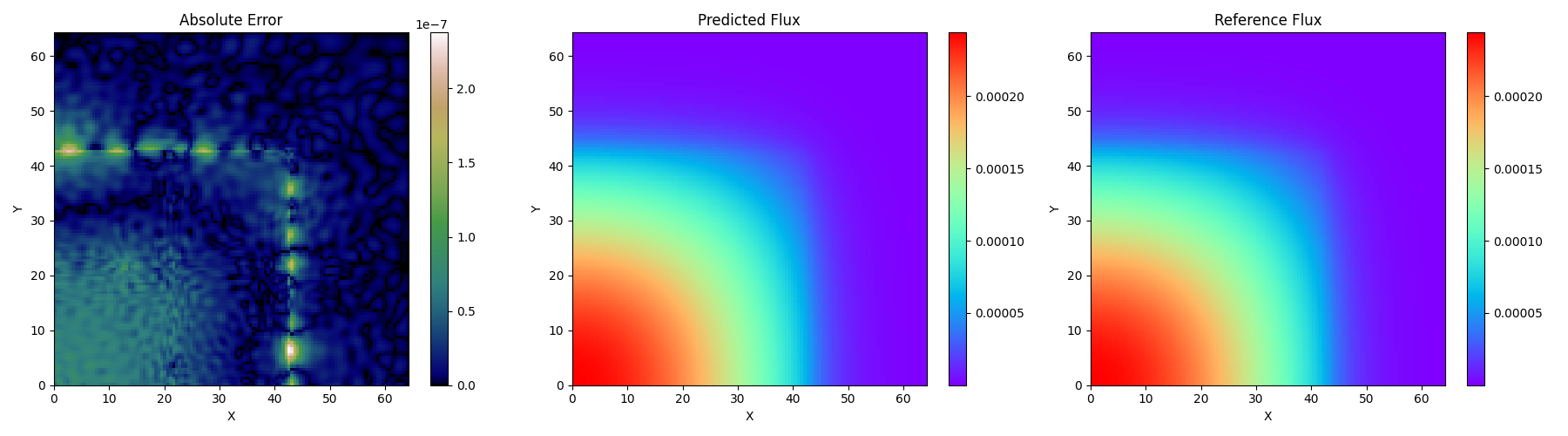}}\\	
	\subfloat[ Thermal group ]{\includegraphics[scale=0.3]{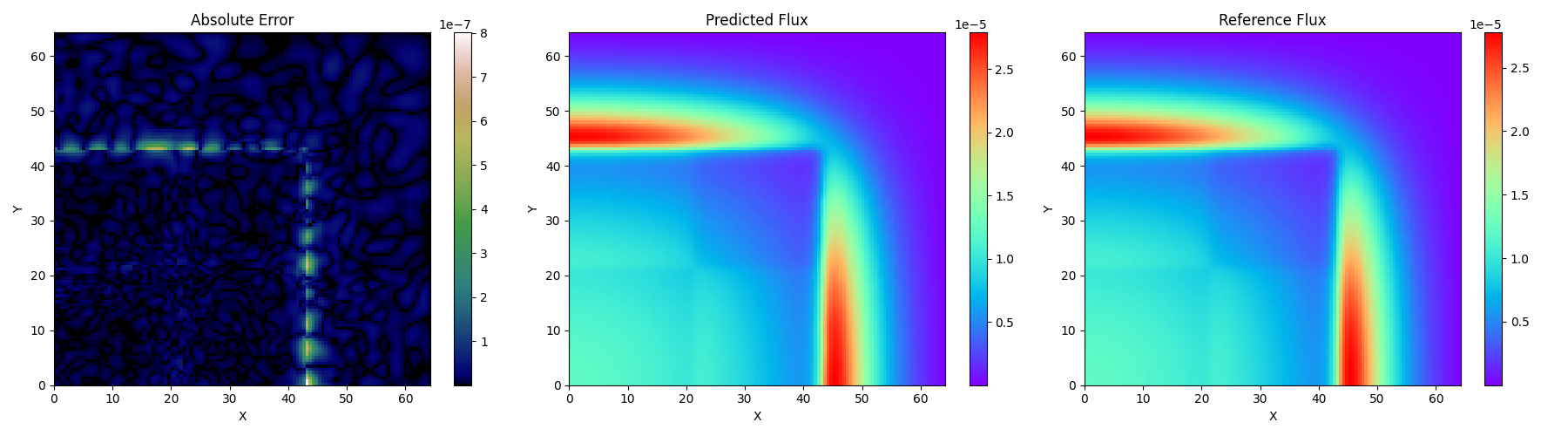}}\\
	
	\caption{ Neutron flux and the deviation from the reference solution for the simplified C5G7 benchmark k-eigenvalue problem with scaled loss function, Sobol sampling method and $Sin$ activation function.}
	\label{fig:C5G7-flux} 
\end{figure}
\subsubsection{The TWIGL-2D benchmark}
\begin{figure}[h]
\centering
\begin{minipage}{0.5\textwidth}
  \centering
  \includegraphics[scale=1.0]{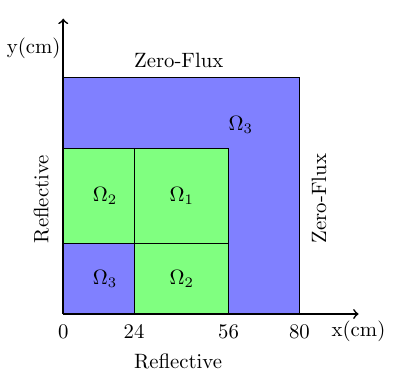} 
  \caption{The TWIGL geometry }
  \label{fig:TWIGL-geo}
\end{minipage}\hfill
\begin{minipage}{0.5\textwidth}
  \centering
   \captionof{table}{Cross sections $(cm^{-1})$ for the TWIGL benchmark} 
  \begin{tabular}{@{}c ccc@{}}
\toprule
 & $\Omega_1$ & $\Omega_2$ & $\Omega_3$ \\
\midrule
$D^1$                     & 1.4  & 1.4  & 1.3 \\
$D^2$                     & 0.4  & 0.4  & 0.5 \\
\addlinespace
$\Sigma_{a}^1$            & 0.01 & 0.01 & 0.008 \\
$\Sigma_{a}^2$            & 0.15 & 0.15 & 0.05 \\
\addlinespace
$\Sigma_{s}^{1\to 2}$        & 0.01 & 0.01 & 0.01 \\
\addlinespace
$\nu\Sigma_{f}^1$         & 0.007 & 0.007 & 0.003 \\
$\nu\Sigma_{f}^2$         & 0.20  & 0.20  & 0.06 \\
 \bottomrule
\end{tabular}

  \label{tab:TWIGL-xs}
\end{minipage}
\end{figure}

First, we would like to consider the TWIGL-2D benchmark, a symmetric seed–blanket problem, whose geometry is shown in Figure \ref{fig:TWIGL-geo} and the corresponding cross sections are presented in Table \ref{tab:TWIGL-xs}.
The calculation of the reference solution is performed on a uniform $120\times120$ grid using the MINOS solver. For this TWIGL benchmark, the reference value of the multiplication factor is $\keffref= 0.91320$ . The number of residual points is $N_r=20480$.

\begin{figure}[H]
	\centering
	\subfloat[$\bar{\Delta} \phi$]{\includegraphics[scale=0.4]{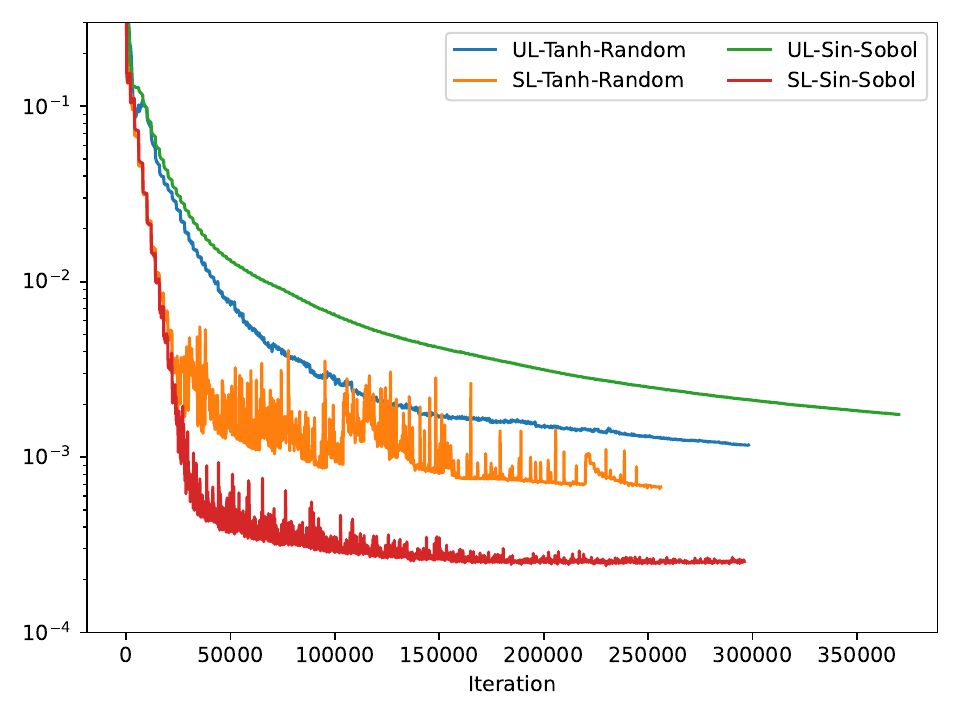}}
    \subfloat[$\bar{\Delta} \pvec$]{\includegraphics[scale=0.4]{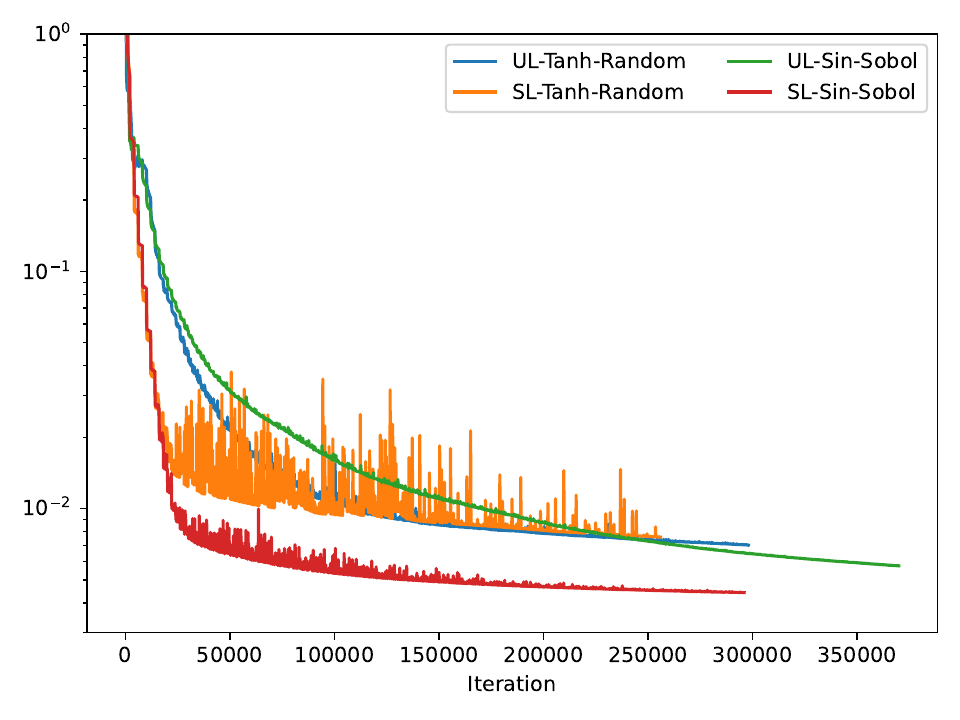}}	
	\caption{ The relative error of the neutron flux and current of the  TWIGL 2D problem.}
    	\label{fig:TWIGL-2D-error} 
\end{figure}

\begin{figure}[H]
	\centering
	\subfloat[Fast group ]{\includegraphics[scale=0.3]{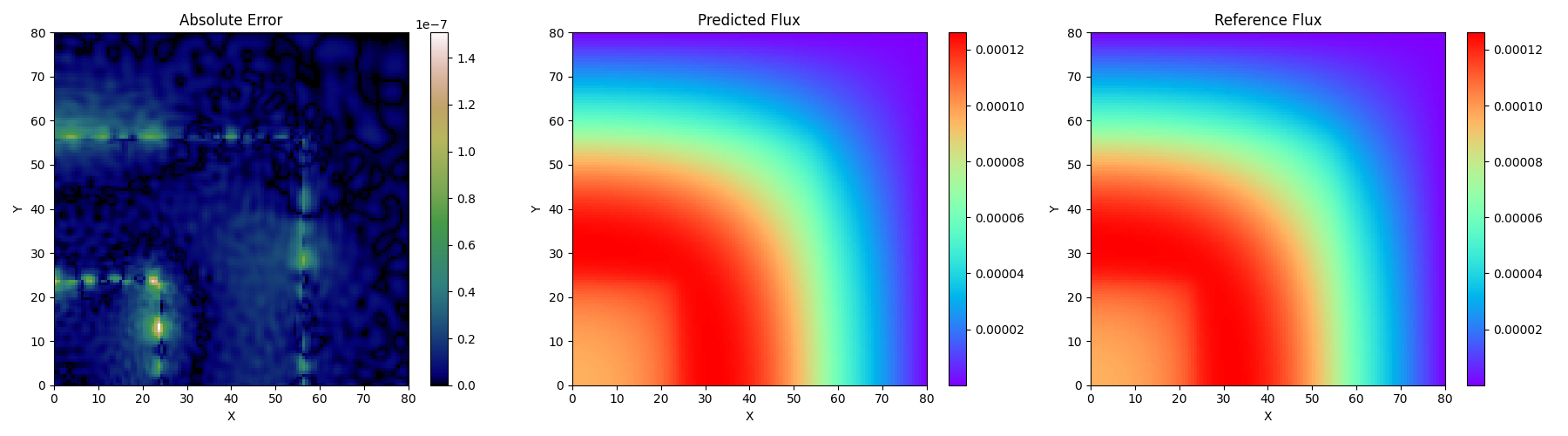}}\\	
	\subfloat[ Thermal group ]{\includegraphics[scale=0.3]{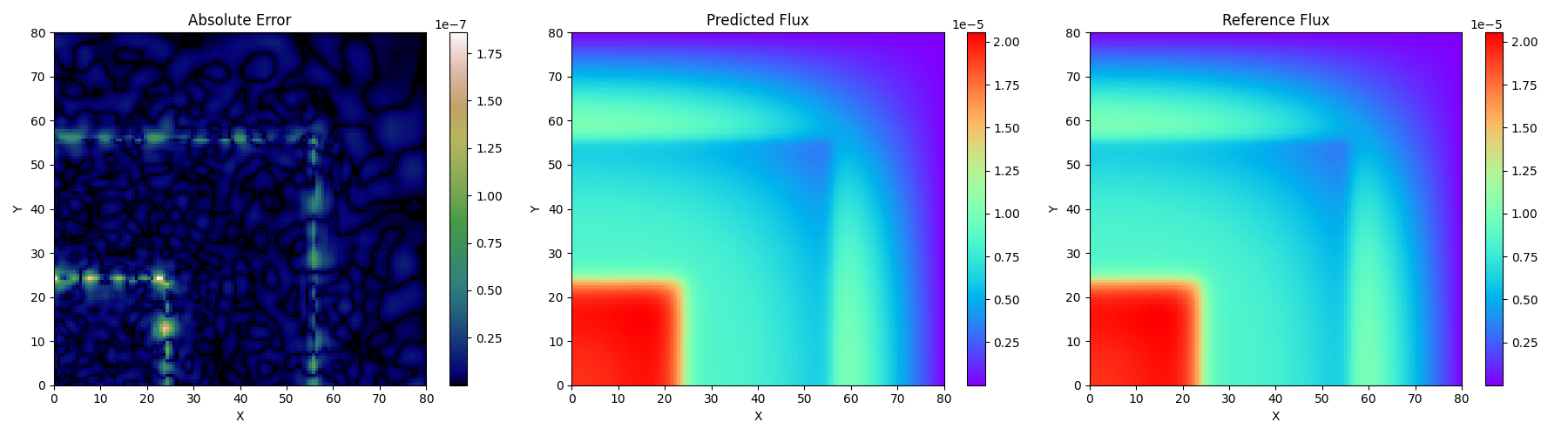}}\\
	
	\caption{ Neutron flux and the deviation from the reference solution for the TWIGL-2D benchmark k-eigenvalue problem with scaled loss function, Sobol sampling method and $Sin$ activation function.}
	\label{fig:TWIGL-2D-flux} 
\end{figure}

\subsubsection{The TWIGL-3D benchmark}
\begin{figure}[h]
\centering
  \includegraphics[scale=0.8]{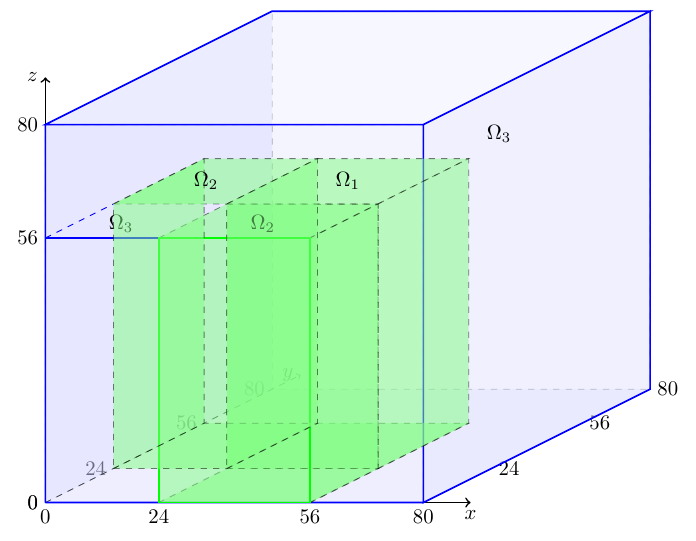} 
  \caption{The TWIGL-3D geometry }
  \label{fig:TWIGL-3D-geo}
  \end{figure}

The TWIGL-3D benchmark is a seed–blanket reactor problem with the geometry given in Figure \ref{fig:TWIGL-3D-geo} and represents a three-dimensional extension of the previous TWIGL-2D benchmark in the z-direction. The same material is used in this test. The TWIGL-3D benchmark uses the same macroscopic cross sections as the TWIGL-2D benchmark (Table \ref{tab:TWIGL-xs}), with the geometry extended in the z-direction.  The reference solution is performed on a uniform $50\times50\times 50$ grid and the reference value of the multiplication factor is $\keffref= 0.88598$. The neural network is trained with $N_r= 200\times 1024$ residual points.

\begin{figure}[H]
	\centering
	\subfloat[$\bar{\Delta} \phi$]{\includegraphics[scale=0.4]{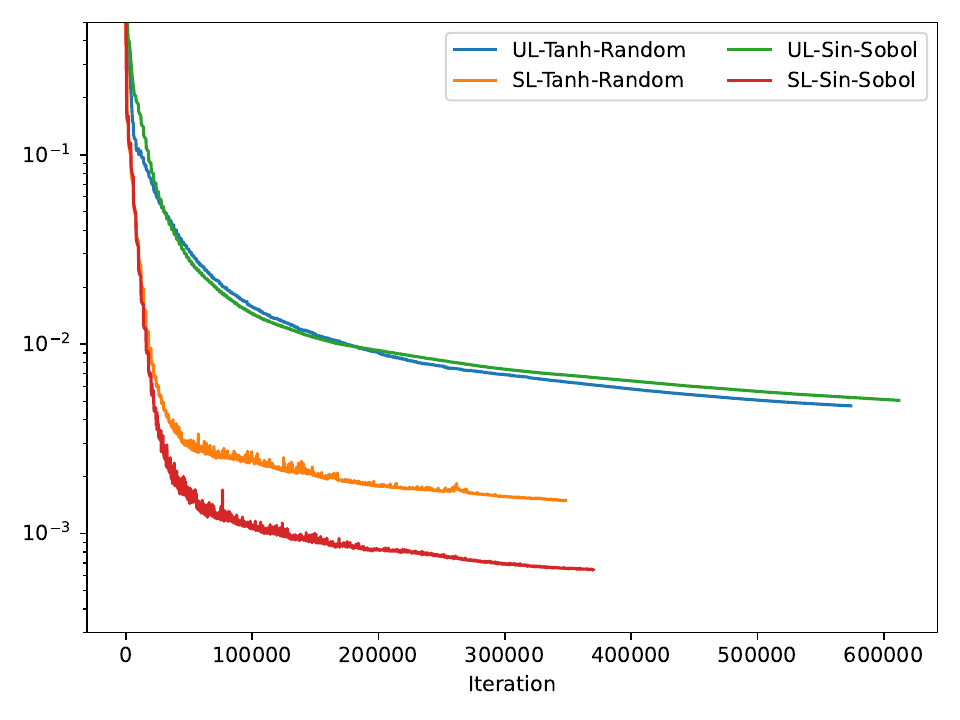}}
    \subfloat[$\bar{\Delta} \pvec$]{\includegraphics[scale=0.4]{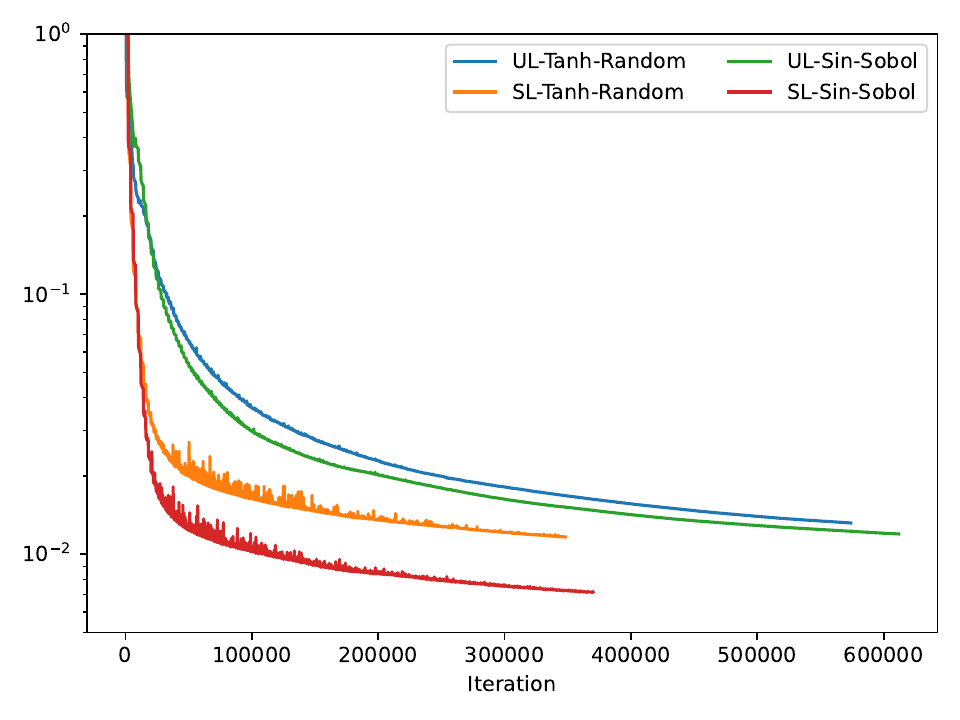}}	
	\caption{ The relative error of the neutron flux and current of the  TWIGL 3D problem.}
    \label{fig:TWIGL-3D-error}
\end{figure}

\begin{figure}[H]
	\centering
	\subfloat[Fast group (a slice in
z-direction at the central index)]{\includegraphics[scale=0.3]{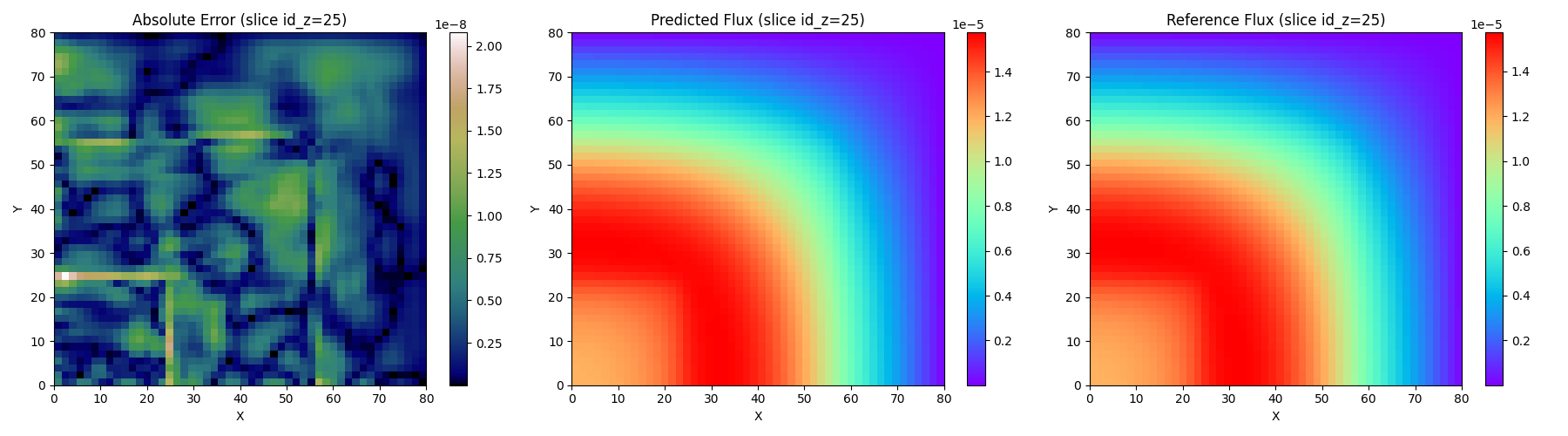}}\\	
	\subfloat[ Thermal group (a slice in z-direction at the central index) ]{\includegraphics[scale=0.3]{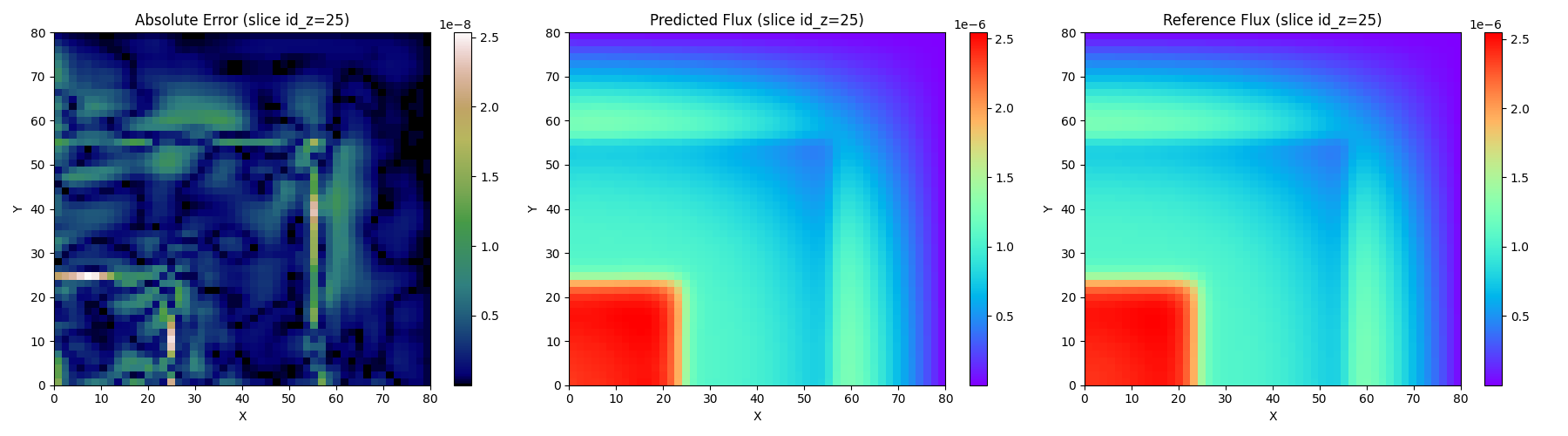}}\\
	
	\caption{ Neutron flux and the deviation from the reference solution for the TWIGL-3D benchmark k-eigenvalue problem with scaled loss function, Sobol sampling method and $Sin$ activation function.}
	\label{fig:TWIGL-3D-flux} 
\end{figure}

Figures \ref{fig:C5G7-error}, \ref{fig:TWIGL-2D-error} and \ref{fig:TWIGL-3D-error} show that the convergence of the scaled loss (SL) is always faster than the standard loss function for all benchmark test cases of the k-eigenvalue problem. The relative error of the neutron flux is always smaller than that of the neutron current. According to Figures \ref{fig:C5G7-flux}, \ref{fig:TWIGL-2D-flux} and \ref{fig:TWIGL-3D-flux}, MF-PINNs with the scaled loss provide accurate approximations for both the fast flux and  thermal flux, with the largest errors concentrated near the material interfaces

Table \ref{tab:errors} summarizes the errors of all test cases for the k-eigenvalue problem. This table shows that the scaled loss (SL) significantly improves the accuracy of $\keff$, neutron flux $\phi$ and neutron current $\pvec$ compared to the unscaled loss (UL). The improvement of the neutron flux ($\phi$) with scaled loss is much more than the neutron current ($\pvec$), especially in the TWIGL 3D case.

The scaled loss function also helps to reduce the number of outer iterations. In addition, the combination of the Sin activation function, Sobol sampling and scaled loss function leads to the lowest error for both the neutron flux and the neutron current across all the tests.

\begin{table}[H]
\centering
\caption{Errors for benchmark test cases }
\label{tab:errors}
\begin{tabular}{cccccccccc}
\toprule
   Test & Method & $N_{\text{outer}}$ & $k_{\text{eff}}$ & $\Delta k_{\text{eff}}$  & $\bar{\Delta}\phi~(\%)$ & $\bar{\Delta}\pvec~(\%)$\\
   \midrule
     \multirow{4}{*}{\rotatebox[origin=c]{90}{SC5G7-2D}} & UL-Tanh-Random & 216 & 0.92701 & 56 & 0.62 &  1.36\\
                                                    & SL-Tanh-Random & 197 & 0.92759 & 2 & 0.09 &0.81  \\
                                                    & UL-Sin-Sobol & 160 & 0.92737 & 20 & 0.15 &  0.57 \\
                                                    & SL-Sin-Sobol & 147 & 0.92760 & 3 & 0.06 & 0.56  \\
\midrule
     \multirow{4}{*}{\rotatebox[origin=c]{90}{TWIGL-2D}} & UL-Tanh-Random & 149 & 0.91288 & 32 & 0.12 & 0.70 \\
                                                    & SL-Tanh-Random & 128 & 0.91316 & 4 & 0.07 & 0.76  \\
                                                    & UL-Sin-Sobol & 185 & 0.91279& 41 & 0.18 & 0.57 \\
                                                    & SL-Sin-Sobol & 148 & 0.91316& 4 & 0.03 & 0.44 \\
   \midrule
  \multirow{4}{*}{\rotatebox[origin=c]{90}{TWIGL-3D}} & UL-Tanh-Random & 287 & 0.88515 & 32 & 0.47 & 1.32 \\
                                                    & SL-Tanh-Random & 174 & 0.88597 & 1 & 0.15 & 1.17  \\
                                                    & UL-Sin-Sobol & 306 & 0.88497 & 101 & 0.51& 1.20 \\
                                                    & SL-Sin-Sobol & 185 & 0.88599& 1 & 0.06& 0.72 \\
\bottomrule
\end{tabular}
\end{table}

\section{Conclusions}
\label{section-conclusion}
This work demonstrates the effectiveness of the Physics Based Loss Scaling (PBLS) method for the mixed formulation of PINNs (MF-PINNs) applied to the neutron diffusion equation. The scaled loss improves the convergence of the MF-PINNs for both the fixed source and the k-eigenvalue problem, compared to the standard loss function. Future works should be dedicated to the investigation of this method in  other 3D test cases, such as the Takeda benchmarks, as well as examining its integration with adaptive sampling approaches \cite{wu2023comprehensive}. Furthermore, extending MF-PINNs to the simplified $P_N$ ($SP_N$) equations offers a promising approach to efficiently solve higher-order approximations of the neutron transport equation.

While MF-PINNs work well for the multigroup neutron diffusion equation, they remain computationally expensive compared to traditional numerical methods and usually require retraining
for each new parameter (cross sections) or configuration. Operator-learning approaches such as Deep Operator Networks (DeepONet) \cite{lu2021learning}, Fourier Neural Operator (FNO) \cite{li2021fourier} and
physics-informed neural operators (PINO) \cite{li2024physics} provide a potential alternative by directly learning the mapping from input parameter space to solution space, enabling a faster approximate solution for  parametric Partial Differential Equations (PDEs).

\section*{Acknowledgement}
APOLLO3\textsuperscript{\textregistered} is a registered trademark of CEA. We gratefully acknowledge EDF and Framatome for their long-term partnership and their support.




 \bibliographystyle{plain} 
 \bibliography{references.bib}







\end{document}